\documentclass{article}

\usepackage{graphicx} 
\usepackage{amsmath,amssymb,amsfonts,amsthm}

\usepackage{soul}
\usepackage{mathtools}
\usepackage{dsfont}
\usepackage{booktabs}
\usepackage{caption}
\usepackage{tablefootnote}
\usepackage{titlesec}
\usepackage{algorithm}
\usepackage{algpseudocode}
\algrenewcommand\Require{\Statex \textbf{Input:}~}
\algrenewcommand\Ensure{\Statex \textbf{Output:}~}
\usepackage{bm}
\usepackage{setspace}
\usepackage{hyperref}
\usepackage{array}
\usepackage{placeins}
\usepackage{xspace}
\usepackage{enumitem}
\usepackage{tikz}
\usetikzlibrary{shapes.geometric, arrows}
\usepackage{xcolor}
\usepackage[sort&compress,numbers]{natbib}

\usepackage{tikz}
\usepackage{url}
\usepackage{threeparttable}
\usepackage{listings}
\usetikzlibrary{calc,positioning,arrows.meta,quotes,decorations.pathreplacing}
\usepackage{pgfplots}
\pgfplotsset{compat=1.17}
\newcommand{\ourmethod}{HALLaR\xspace}
\newcommand{\gpuourmethod}{cuHALLaR\xspace}
\newcommand{\culorads}{cuLoRADS\xspace}

\DeclareMathOperator*{\argmin}{arg\,min}
\algrenewcommand\Require{\State \textbf{Input: }}
\algrenewcommand\Ensure{\State \textbf{Output: }}

\usepackage{amsmath}
\usepackage{amsfonts}
\usepackage{latexsym}
\usepackage{amssymb}
\usepackage{stmaryrd}
\usepackage{mathtools}

\newtheorem{thm}{Theorem}[section]

\newtheorem{remark}[thm]{Remark}

\newcommand{\beq}{\begin{equation}}
\newcommand{\eeq}{\end{equation}}
\newcommand{\beqa}{\begin{eqnarray}}
\newcommand{\eeqa}{\end{eqnarray}}
\newcommand{\beqas}{\begin{eqnarray*}}
\newcommand{\eeqas}{\end{eqnarray*}}
\newcommand{\ei}{\end{itemize}}

\def\defi{\vcentcolon=}

\newcommand{\R}{\mathbb{R}}

\renewcommand{\SS}{\mathbb{S}}

\newcommand{\cA}{{\cal A}}

\newcommand{\lam}{{\lambda}}

\newcommand{\inner}[2]{\langle #1,#2\rangle}

\newcommand{\ignore}[1]{}

\newtheorem{innercustomgeneric}{\customgenericname}
\providecommand{\customgenericname}{}
\newcommand{\newcustomtheorem}[2]{%
  \newenvironment{#1}[1]
  {%
   \renewcommand\customgenericname{#2}%
   \renewcommand\theinnercustomgeneric{##1}%
   \innercustomgeneric
  }
  {\endinnercustomgeneric}
}
\newcustomtheorem{blackbox}{Blackbox}

\tikzstyle{startstop} = [
  rectangle, rounded corners,
  minimum width=4cm, minimum height=1cm,
  text width=6cm, align=justify,
  text centered, draw=blue!80, fill=blue!10
]
\tikzstyle{process} = [rectangle, minimum width=3cm, minimum height=1cm, text centered, draw=purple!80, fill=purple!10]
\tikzstyle{arrow} = [thick,->,>=stealth]
\tikzstyle{bigprocess} = [
  rectangle, rounded corners,
  minimum width=6cm,
  minimum height=1.5cm,
  text width=6cm,
  align=justify,
  draw=purple!80,
  fill=purple!10,
  font=\small          
]
\tikzstyle{arrow} = [thick,->,>=stealth]

\allowdisplaybreaks

\title{An Accelerated Hybrid Low-Rank Augmented Lagrangian Method for Large-Scale Semidefinite Programming on GPU}

\author{
{Jacob M. Aguirre,}
\thanks{H. M. Stewart School of Industrial and Systems Engineering, Georgia Tech, Atlanta, GA, 30332-0205. (Email: {\tt aguirre@gatech.edu})}
{\quad Diego Cifuentes,}
\thanks{H. M. Stewart School of Industrial and Systems Engineering, Georgia Tech, Atlanta, GA, 30332-0205. (Email: {\tt dfc3@gatech.edu})}
{\quad Vincent Guigues}
\thanks{School of Applied Mathematics, FGV,
Praia de Botafogo, Rio de Janeiro, Brazil. (Email:{\tt vincent.guigues@fgv.br})}
{Renato D.C. Monteiro,}
\thanks{H. M. Stewart School of Industrial and Systems Engineering, Georgia Tech, Atlanta, GA, 30332-0205. (Email: {\tt rm88@gatech.edu})}
{\quad Victor Hugo Nascimento,}
\thanks{EMAp FGV.}
{\quad Arnesh Sujanani.}
\thanks{Department of Combinatorics and Optimization, University of Waterloo, Waterloo, ON, N2L 3G1.(Email: {\tt a3sujana@uwaterloo.edu}) \\~\\~\\~\\ Authors are listed in alphabetical order}
}

\date{May 19, 2025 (second version: October 14th, 2025)}

\title{
\gpuourmethod: A GPU Accelerated Low-Rank Augmented Lagrangian Method for Large-Scale Semidefinite Programming
\thanks{\textbf{Funding}: 
Jacob M. Aguirre is supported by the National Science Foundation Graduate Research Fellowship under Grant No. DGE-2039655.
Diego Cifuentes is supported by U.S. Office of Naval Research, N00014-23-1-2631.
Renato D.C. Monteiro is supported by AFOSR Grant FA9550-25-1-0131. 
}}

\begin{document}
\author{Jacob M. Aguirre\thanks{H. M. Stewart School of Industrial and Systems Engineering, Georgia Tech, Atlanta, GA, 30332-0205. \protect\protect\href{mailto:aguirre@gatech.edu}{aguirre@gatech.edu}}\hspace*{0.5em}
\and Diego Cifuentes\thanks{H. M. Stewart School of Industrial and Systems Engineering, Georgia Tech, Atlanta, GA, 30332-0205.  
\protect\protect\href{mailto:dfc3@gatech.edu}{dfc3@gatech.edu}}
\and Vincent Guigues\thanks{School of Applied Mathematics, FGV,
Praia de Botafogo, Rio de Janeiro, Brazil. 
\protect\protect\href{mailto:vincent.guigues@fgv.br}{vincent.guigues@fgv.br}}
\and Renato D.C. Monteiro\thanks{School of Industrial and Systems Engineering, Georgia Institute of
Technology, Atlanta, GA, 30332-0205. 
\protect\protect\href{mailto:monteiro@isye.gatech.edu}{monteiro@isye.gatech.edu}}
\and Victor Hugo Nascimento\thanks{School of Applied Mathematics, FGV,
Praia de Botafogo, Rio de Janeiro, Brazil.
\protect\protect\href{mailto:nascimento.victor.1@fgv.edu.br}{nascimento.victor.1@fgv.edu.br}}
\and Arnesh Sujanani\thanks{Department of Combinatorics and Optimization, University of Waterloo, Waterloo, ON, N2L 3G1.
\protect\protect\href{mailto:a3sujana@uwaterloo.ca}{a3sujana@uwaterloo.ca}\\~\\ Authors are listed in alphabetical order.}
}

\setcitestyle{square}
\maketitle

\begin{abstract}
This paper presents an SDP solver, \gpuourmethod, which is
a GPU-accelerated implementation of the hybrid low-rank augmented Lagrangian method, \ourmethod, proposed earlier by three of the authors. The proposed Julia-based implementation is designed to exploit massive parallelism by performing all core operator evaluations—linear maps, their adjoints, and gradients—entirely on the GPU. Extensive numerical experiments are conducted on three problem classes: matrix completion, maximum stable set, and phase retrieval. The results demonstrate substantial performance gains over both optimized CPU implementations and existing GPU-based solvers. On the largest instances, \gpuourmethod achieves speedups of up to 165x, solving a matrix completion SDP with a matrix variable of size $2 \times 10^6 \times 2 \times 10^6$ and over $2.6 \times 10^8$ constraints to a relative precision of $10^{-5}$ in 53 seconds. Similarly, speedups of up to 135x are observed for maximum stable set problems on graphs with over one million vertices. These results establish \gpuourmethod as a state-of-the-art solver for extremely large-scale SDPs.
\\


\end{abstract}

\section{Introduction}\label{Introduction}

Let $\mathbb S^n$ be the set of
$n\times n$ symmetric matrices.
The notation $A \succeq B$ means
that $A-B$ is positive semidefinite.
We are interested in solving the primal-dual pair of semidefinite programs (SDPs)
\begin{gather*}
\label{eq:sdp-primal}\tag{P}
P_{*}: = \quad \min_{X} \quad \{C \bullet X
\quad : \quad
\mathcal A(X)=b ,\quad
X  \in \Delta^n_{\tau} \}
\end{gather*}
and
\begin{gather*}
\label{eq:sdp-dual}\tag{D}
D_{*}: = \quad \max_{p \in \mathbb{R}^m, \theta \in \mathbb{R}_+}   \quad \{ -b^{\top}p-\tau^2\theta
\quad : \quad
S \defi C+\mathcal A^{*}(p)+\theta I \succeq 0\}
\end{gather*}
where $\bullet$ is the Frobenius inner product, $b \in \mathbb{R}^m$, $C\in \mathbb S^n$, and $\mathcal A: \mathbb S^n \to \mathbb{R}^m$ and $\mathcal A^*: \mathbb{R}^m \to \mathbb S^n$ are linear maps such that 
\begin{align}\label{Definition of Linear Operators}
     \mathcal A(X)= \begin{bmatrix}
           A_1\bullet X \\
           A_2\bullet X \\
           \vdots \\
           A_m \bullet X
         \end{bmatrix}, \quad \mathcal A^{*}(p)=\sum_{\ell=1}^{m}p_\ell A_\ell, \quad \text{where $A_{\ell} \in \mathbb S^n$ for $\ell=1,\ldots m$} .
        \end{align}
We also define $\Delta^n_{\tau}$ to be the spectraplex
\begin{equation}\label{Delta Definition}
\Delta^n_{\tau} \defi \{ X \in \mathbb S^n : \mathrm{Tr}(X) \leq \tau^2, X \succeq 0 \}.
\end{equation}
Semidefinite programming arises across a wide spectrum of applications, notably in machine learning, computational sciences, and various engineering disciplines. Despite its versatility, efficiently solving large-scale SDPs remains computationally demanding. Interior-point methods, although widely employed, are often constrained by computational complexity as they often rely on solving a large linear system using a Newton-type method, which involves an expensive inversion of a large matrix
\cite{borchers1999csdp, Overton, ipmWolk, Mon97, ToddToh}. These computational bottlenecks have sparked interest in developing efficient first-order methods for solving large-scale problems. Unlike interior-point methods, first-order methods scale well with the dimension of the problem instance and predominantly utilize sparse matrix-vector multiplications to perform function and gradient evaluations.
For an extensive survey and analysis of recent developments in using first-order methods for solving SDPs, we refer the reader to the comprehensive literature available in \cite{Alfakih, deng2022enhanced, DingGrimmer, ding2021optimal, douik2018low, erdogdu2022convergence, garstka2021cosmo, huang2017solving, DNNPeaceman, li2021strictly, Kang, o2016conic, pham2023scalable, renegar2019accelerated, shinde2021memory, wen2010alternating, yang2015sdpnal+, yurtsever2019conditional, yurtsever2021scalable, zhao2010newton, zheng2017fast}.

Due to GPU's great capabilities in parallelizing matrix-vector multiplications, they have been shown to work extremely well in speeding up first-order methods in practice \cite{chen2024cuclarabel, han2024accelerating, lin2025pdcs, lu2023cupdlpC, lu2023practical, lu2023cupdlp, schubiger2020gpu}. In this paper, we introduce an enhanced Julia-based GPU implementation, referred to as 
\gpuourmethod,
of the first-order \ourmethod method developed in \cite{monteiro2024low} for solving huge-scale SDPs. Our work reinforces the viability of using GPU architectures for efficient parallelization in semidefinite programming solvers. Moreover, our experiments show potential for excellent speedups when using GPUs for optimization solvers. 

\subsection{Related literature}

We divide our discussion into four parts: first, the emerging literature on utilizing GPUs for mathematical optimization, second, the challenges in solving large SDPs, third, recent techniques in solving large SDPs using low-rank approaches, and fourth, a comparison between \gpuourmethod, a recent GPU SDP solver \culorads developed in \cite{han2024accelerating}, and the CPU version of \ourmethod. 

\paragraph{GPUs for Mathematical Programming}
GPUs offer massive parallelism that accelerates mathematical programming through parallel linear algebra operations. Recent GPU-based solvers for linear programming and quadratic programming \cite{lu2023cupdlp,lu2023cupdlpC, lu2023practical} have achieved performance which outperforms traditional CPU-based methods. First-order methods map well to GPU architectures, unlike simplex or interior-point methods which rely on difficult-to-parallelize factorizations. For linear programming, cuPDLP.jl incorporates GPU-accelerated sparse matrix operations to match commercial solver performance. Quadratic programming solvers like cuOSQP \cite{schubiger2020gpu} and CuClarabel \cite{chen2024cuclarabel} similarly exploit GPU parallelism for significant speedups. 

Most recently, \cite{lin2025pdcs} proposed a new GPU primal-dual conic programming solver, namely, PDCS, which utilizes the primal-dual hybrid gradient (PDHG) method for solving several classes of cone optimization problems, including semidefinite programs; however, the authors of the paper remark that due to PDHG requiring Euclidean projections, their method for solving SDP problems is rather inhibited and slow. Thus, their code does not allow users to input SDP data, and so we cannot compare against PDCS. 

\paragraph{Challenges in large-scale semidefinite programming}
Solving large-scale SDPs in practice presents greater challenges compared to solving large-scale linear programs and quadratic programs due to the reliance of many SDP solvers on full 
eigendecomposition, inversion, and storage of large dense matrices. To circumvent several of these issues, pre-processing techniques such as facial reduction, chordal conversion, and symmetry reduction have been proposed which can potentially convert the original SDP \eqref{eq:sdp-primal} into a smaller one or one with a more favorable structure 
\cite{vandenberghe2015chordal, kojima, Wolk, zhang2021sparse, zheng2020chordal, zhang2024complexity, hu2023facial, borwein1981regularizing, drusvyatskiy2017many, tang2024exploring, permenter2018partial}. 
Interior-point methods are then often applied to the smaller reformulated  SDP. Another popular approach  that is commonly
used in practice to solve large-scale SDPs and which avoids storing and inverting
large dense matrices is the low-rank approach first proposed by Burer and 
Monteiro \cite{burer2003nonlinear, burer2005local}. The next paragraph 
describes low-rank first-order methods for solving SDPs in more detail.

\paragraph{Low-rank first-order methods for SDPs}
The low-rank approach proposed by Burer and Monteiro \cite{burer2003nonlinear, burer2005local} converts SDP \eqref{eq:sdp-primal} into a nonconvex quadratic program with quadratic constraints (QCQP) through the transformation $X=UU^{\top}$, where $U \in \mathbb R^{n \times r}$ and $r\ll n$. The main advantage of this approach is that the QCQP has $n\times r$ variables, while \eqref{eq:sdp-primal} has $n(n+1)/2$ variables. Also, the PSD constraint $X\succeq 0$ has been eliminated through the low-rank transformation so methods do not need to perform eigendecompositions of dense matrices to approximately solve the QCQP. Burer and Monteiro proposed an efficient first-order limited-memory BFGS augmented Lagrangian (AL) method, which only relied on matrix-vector multiplications and storage of $n\times r$ matrices, to find an approximate local minimizer of the QCQP. Recent landscape results have shown that if the factorization parameter $r$ is taken large enough, then finding local minimizers of the QCQP is actually equivalent to finding global minimizers of the SDP \eqref{eq:sdp-primal} \cite{cifuentes2021burer, cifuentes2022polynomial, boumal2016non, boumal2020deterministic, bhojanapalli2018smoothed, pumir2018smoothed, waldspurger2020rank}. For a more extensive list of papers which characterize the optimization landscape of nonconvex formulations obtained through Burer-Monteiro factorizations, see \cite{ouyang2025burer, endor2024benign, mcrae2025benign, ling2025local, o2022burer, ding2025squared, mcrae2024benign, levin2025effect}. For recent advances and developments in accelerating low-rank first-order methods for solving SDPs, we refer the reader to the extensive literature available in \cite{ erdogdu2022convergence, han2024low, han2024accelerating, monteiro2024low, TangToh, TangToh2, TangToh3, wang2023solving, RiemannianSemiSmooth, yurtsever2015scalable, wang2025accelerated, yurtsever2019conditional, yurtsever2021scalable}.

We now give a brief comparison between \ourmethod, 
which was first developed in 
\cite{monteiro2024low} and the LoRADS method 
developed in \cite{han2024low}. \ourmethod 
utilizes an inexact AL method to solve \eqref{eq:sdp-primal}. It solves a sequence of AL subproblems of the form $\argmin_{X\in \Delta^{n}_\tau}\
    C \bullet X + p^\top (\mathcal{A}(X) - b) + \beta \| \mathcal{A}X - b\|^2/2$ by first restricting them to the space of 
    matrices of the form $X=UU^{\top}$, where the number of columns $r$ of $U$ is significantly smaller than $n$,
    and then applying a 
    nonconvex solver to the reformulation. It uses a low-rank Frank-Wolfe step to escape 
    from a possible spurious stationary point obtained by the nonconvex solver, a step which generally 
    adds an additional column to $U$.
    
    LoRADS, on the other hand, adopts a two-stage approach. In their 
    first stage, they utilize the low-rank 
    limited-memory BFGS AL method of Burer and Monteiro as a warm-start strategy for their second stage. 
    The second stage of
    LoRADS also utilizes an inexact AL method applied to \eqref{eq:sdp-primal} but the major difference being the way 
    the AL subproblems are reformulated. In contrast to the $X=UU^{\top}$ change of variables used by \ourmethod, LoRADS
    reformulates them 
by replacing
    $X$ by 
    $UV^{\top}$ and adding the constraint 
    $U=V$.
    It then approximately solves these nonconvex reformulations of the AL subproblems by an alternating minimization penalty-type approach which penalizes $U=V$ constraint and alternately minimizes with respect to the blocks
    $U$ and $V$ using the conjugate gradient method.
    Similar to \ourmethod,
    the inexact stationary solutions of the nonconvex AL  reformulations may not be near global solutions of the AL subproblems as they are not necessarily equivalent.
   In contrast to \ourmethod, LoRADS 
   uses a heuristic that increases the number of columns of $U$ and $V$
    to attempt to escape from the these spurious solutions.

\paragraph{\gpuourmethod versus cuLoRADS and \ourmethod}
Like cuLoRads, which is a GPU accelerated version of LoRads developed in \cite{han2024accelerating},  \gpuourmethod is also a GPU accelerated version of \ourmethod. 
Computational results reported in this paper show that 
\gpuourmethod is 2 to 25 times faster than cuLoRads
on large-scale instances of some important SDP classes. \gpuourmethod can efficiently solve massive problems to $10^{-5}$ relative precision in just a few minutes. For example, \gpuourmethod is able to solve a matrix completion SDP instance with matrix variable of size 8 million, and  approximately 300 million constraints, in just 142 seconds.

\gpuourmethod also achieves massive speedup compared to its CPU counterpart, \ourmethod. \gpuourmethod achieves speedups of 30-140x on large matrix completion SDP instances, up to 135x speedup on large maximum stable set SDP instances, and 15-47x speedup on phase retrieval SDP instances. These numerical results show that \gpuourmethod is a very promising GPU-based method for solving extremely large SDP instances.

\subsection{Notation}
Let $\mathbb{R}^n$ be the space of $n$ dimensional vectors,
$\mathbb{R}^{n\times r}$ the space of $n \times r$ matrices, and $\mathbb{S}^n$ the space of symmetric $n \times n$ matrices. Let $\mathbb{R}_{++}^n$ (resp., $\mathbb{R}^n_{+}$) denote the convex cone in $\mathbb{R}^n$ of vectors with positive (resp., nonnegative) entries
and  $\inner{\cdot}{\cdot}$ (resp. $\|\cdot\|$) be the Euclidean inner product (resp. norm) on~$\mathbb{R}^n$. 
Given matrices $A_1,\ldots,A_m \in \mathbb{S}^n$, let 
$\mathcal{A}: \mathbb{S}^n \rightarrow \mathbb{R}^m$ denote the operator defined as
$[\mathcal{A}(X)]_i=A_{i} \bullet X$ for every $i=1,\ldots,m$.
The adjoint operator $\mathcal{A}^*:\mathbb{R}^m \rightarrow \mathbb{S}^n$ of 
$\mathcal{A}$ is given by 
$\mathcal{A}^*(\lambda)=\sum_{i=1}^m \lambda_i A_i$. The minimum eigenvalue of a matrix $Q \in \mathbb S^n$ is denoted by $\lambda_{\min}(Q)$, and $v_{\min}(Q)$ denotes a corresponding eigenvector of unit norm.

For a given $\epsilon\geq 0$, the $\epsilon$-normal cone of a closed convex set $C$ at $z\in C$, denoted by  $N^{\epsilon}_C(z)$, is 
$$N^{\epsilon}_C(z):=\{\xi \in \mathbb R^n: \inner{\xi}{u-z}\leq \epsilon,\quad \forall u\in C\}.$$
The normal cone of a closed convex set $C$ at $z\in C$ is denoted by  $N_C(z)=N^0_C(z)$. Finally, we define the Frobenius ball of radius $r$ in $\mathbb R^{n \times s}$ space to be
\begin{equation}\label{dimensional ball}
B_r^s:=\{U \in \mathbb R^{n \times s}:\|U\|_{F}\leq r\}.
\end{equation}
We denote the $i$-th row of a matrix $U$ by $U_{i,:}$.

\section{Overview of \ourmethod}\label{Overview}
This section provides a minimal review of \ourmethod, the hybrid low-rank augmented Lagrangian method that was first developed in \cite{monteiro2024low}. 

As mentioned in the introduction, \ourmethod is an inexact augmented Lagrangian method that approximately solves the pair of semidefinite programs \eqref{eq:sdp-primal} and \eqref{eq:sdp-dual}. Given a tolerance pair $(\epsilon_p,\epsilon_{pd}) \in \mathbb{R}^2_{++}$, the method aims to find an $(\epsilon_p,\epsilon_{pd})$-optimal solution, i.e., a triple $(\bar X, \bar p, \bar \theta) \in \Delta^{n}_\tau \times \mathbb{R}^m \times \mathbb{R}_{+}$ that satisfies
\begin{equation}\label{SDP criteria}
\|\mathcal A\bar X-b\|\leq \epsilon_p, \quad |C \bullet \bar X + b^\top \bar p + \tau^2\bar\theta|\leq \epsilon_{pd}, \quad C+\mathcal A^{*}\bar p+\bar \theta I \succeq 0, \quad
\bar \theta \geq 0.
\end{equation}
When the tolerances are zero, these relations reduce to the standard optimality conditions for the primal-dual pair. Note that the final two inequalities in \eqref{SDP criteria} enforce dual feasibility, a property that \ourmethod maintains at every iteration. Being an inexact AL method, \ourmethod generates sequences $\{X_t\}$ and $\{p_t\}$ according to the following updates
\begin{subequations}
\begin{align}
    \label{eq:AL}
    X_{t} \quad &\approx \quad
    \argmin_X\quad
    \{\mathcal{L}_{\beta_t}(X;p_{t-1}) \;:\; X \in \Delta^n_{\tau}\},
    \\
    p_{t} \quad &= \quad
    p_{t-1} + \beta_t(\mathcal{A}(X_{t}) - b)\label{lm update}
\end{align}
\end{subequations}
where
\begin{align}\label{AL function}
\mathcal{L}_\beta(X;p)
    \quad\defi\quad
    C \bullet X + p^\top (\mathcal{A}(X) - b) + \frac{\beta}{2} \| \mathcal{A}X - b\|^2
\end{align}
is the AL function.
The key part of \ourmethod is an efficient method, called the hybrid low-rank method (HLR), for finding a near-optimal
solution $X_{t}$ of the AL subproblem~\eqref{eq:AL}. HLR never forms $X_{t}$ but outputs a low-rank factor $U_t \in B_{\tau}^{s_t}$ such that $X_t=U_tU_t^{\top}$, where $B_{\tau}^{s_t}$ is as in \eqref{dimensional ball} and hopefully $s_t\ll n$. The next paragraph briefly describes HLR in more detail for solving \eqref{eq:AL}.

Given a pair $(s, \tilde Y) \in \mathbb Z_{+} \times B_{\tau}^s  $, which is initially set to $(s,\tilde Y)=(s_{t-1},U_{t-1})$,
a general iteration of HLR performs the following steps:
i) it computes
a suitable stationary point $Y$ of the nonconvex reformulation of \eqref{eq:AL}
obtained through the change of variable $X=UU^{\top}$, namely:
\begin{align}\label{eq:Lr}
  \tag{$\mathrm{L}_r$}
  \min_U\quad
  \left \{ 
    g_t(U) = \mathcal{L}_{\beta_t}(UU^\top;p_{t-1}) \quad:\quad
    \|U\|_F \le \tau, \quad U \in \mathbb{R}^{n\times s} \right\};
\end{align}
(ii) it verifies if 
the point 
$X=YY^{\top}$ 
corresponds to a near-global minimizer of the original convex AL subproblem \eqref{eq:AL}. This check is performed by computing the gradient $G$ of the AL function at $X=YY^\top$ and using its minimum eigenpair to compute a Frank-Wolfe optimality gap. If this gap is below a prescribed tolerance, the point is accepted as a suitable solution of \eqref{eq:Lr}, and HLR stops, setting $(s_t,U_t)=(s,Y)$;
(iii) otherwise, if the gap is not small then this indicates that $Y$ is a spurious stationary point of \eqref{eq:Lr}. The minimum eigenvector computed in the previous step furnishes a descent direction for the convex problem \eqref{eq:AL}. Using this direction, a Frank-Wolfe step is performed in the $X$-space (but implemented in the $U$-space) to escape the spurious point $Y$. This step results in a rank-one update that produces a new iterate pair $(s,\tilde Y)$.





We now give more details about the above steps.
First, we consider step i).
This step is implemented with the aid of a generic nonlinear solver which, started from 
$\tilde Y$,  computes a
$(\rho;\tilde Y)$-stationary solution of \eqref{eq:Lr}, i.e., a triple $(Y;R,\rho)$ such that

\begin{equation}\label{Stationary}
    R \in \nabla g_t(Y)+N_{B^s_{\tau}}(Y), \quad \|R\|_{F}\leq \rho, \quad g_t(Y)\leq g_t(\tilde Y).
\end{equation}
It is easy to verify that 
\eqref{Stationary}, in terms of the SDP data, is equivalent to
\[
 R \in 2[C+\mathcal A^{*}(p_{t-1}+\beta_t(\mathcal A(YY^{\top})-b))]Y+N_{B^s_{\tau}}(Y), \quad \|R\|_{F}\leq \rho.
\]
Appendix~\ref{ADAP-AIPP} describes a method that is able to carry out step i), namely,
ADAP-AIPP, whose exact formulation and analysis can be found in \cite{monteiro2024low, SujananiMonteiro}, and also in earlier works (see e.g. \cite{Aaronetal2017, CatalystNC, WJRproxmet1, WJRComputQPAIPP}) in other related forms.

The details of steps ii) and iii) implemented in the $X$-space are quite standard.
The details
of its implementation in the $U$-space are described in the formal description of the algorithm given below.
The only point worth observing at this stage is that, due to the rank-one update nature of the Frank-Wolfe method, at the end of step iii), the low-rank parameter $s$ is set to either
$s=1$ (unlikely case) or $s=s+1$ (usual case).
Hence, in the usual case, the number of columns of $\tilde Y$ increases by one.

We now formally state \ourmethod below.
\begin{algorithm}[H]
\setstretch{0.8}
\caption{\normalsize \ourmethod.}
\begin{algorithmic}[1]
\normalsize \Require Let initial points $(U_0, p_0) \in B^{s_0}_{\tau} \times \mathbb{R}^m$, 
tolerance pair $(\epsilon_{\mathrm{p}}, \epsilon_{\mathrm{pd}}) \in \mathbb{R}^2_{++}$, let $\{\beta_{t}\}_{t\geq 1}$ be a sequence of positive reals, 
and let $\left\{\epsilon_t \right\}_{t\geq 1}$ 
be a decreasing 
sequence of positive reals converging to $0$.
\Ensure $(\bar X, \bar p, \bar \theta) \in \Delta^{n}_\tau \times \mathbb{R}^m \times \mathbb{R}_{+}$, an $(\epsilon_{\mathrm{p}},\epsilon_{\mathrm{pd}})$-solution of the pair of SDPs \eqref{eq:sdp-primal} and \eqref{eq:sdp-dual}.
\State set $t \gets 1$;
\State set $\tilde Y=U_{t-1}$, $s=s_{t-1}$, $Y=0_{n\times s}$, $G=0_{n\times n}$, and $\theta=\infty$;
\While{$G \bullet (YY^\top) + \tau^2\theta >  \epsilon_t$ }
\Comment{(\textbf{HLR method)}}
\State call a nonconvex solver with initial point $\tilde Y$, tolerance $\epsilon_t$, and function $\mathcal{L}_{\beta_t}(UU^\top;p_{t-1})+\delta_{B_{\tau}^s}(U)$
to find a point $Y \in B^{s}_{\tau}$ that is an $\epsilon_t$-approximate stationary solution (according to the criterion in \eqref{Stationary}) of
\begin{align}\label{Nonconvex subproblem}
  \min_U\quad
  \left \{ 
\mathcal{L}_{\beta_t}(UU^\top;p_{t-1}) \quad:\quad
    \|U\|_F \le \tau, \quad
  U \in \mathbb{R}^{n\times s} \right\};
\end{align}
\State compute
\begin{equation}\label{gradient Def Algo 2}
G := \nabla \left[ \mathcal L_{\beta_t}\left(\,\cdot\,;p_{t-1}\right)\right]\left(YY^{\top}\right)
\in \mathbb S^{n}
\end{equation}
and a minimum eigenpair of $G$, i.e., $(\lambda_{\min}(G), v_{\min}(G)) \in \mathbb{R} \times \mathbb{R}^n$, and set
\begin{equation}\label{FW subproblem eq Algo 2}
\theta := \max\{-\lambda_{\min}(G) , 0\}, \qquad
y := \begin{cases}
     v_{\min}(G) & \text{ if $\theta>0$}\\
        0 & \text{ otherwise};
        \end{cases}
\end{equation}
\If{$G \bullet (YY^\top) + \tau^2\theta \leq \epsilon_t$}
    \State \textbf{break while loop} and go to step 12 below;
\EndIf

\State compute 
\begin{align}\label{FW stepsize}
  \alpha = \min \left \{ \frac{G \bullet (YY^\top) + \tau^2 \theta}{\beta_t\|\mathcal{A}(Y Y^\top) - \mathcal{A}(\tau^2 y y^\top)\|^2}, \; 1\right\}
\end{align}
and set 
\begin{equation}\label{tilde Y Algo 2}
(\tilde Y,s)=\begin{cases}
(\tau y,1) & \text{ if $\alpha=1$}\\
\left(\left[\sqrt{1-\alpha} \, Y , \sqrt{\alpha} \, \tau \, {y}\right],s+1\right)& \text{ otherwise};
\end{cases}
\end{equation}
\EndWhile
\State set $(U_t,\theta_t,s_t)=(Y,\theta,s)$;
\State set
\begin{equation}\label{dual update SDP Algo 2}
p_t=p_{t-1}+\beta_t (\mathcal A(U_t U_t^\top)-b);
\end{equation}
\State if 
\begin{equation}\label{termination}
\|\mathcal{A}(U_t U_t^\top) - b\| \leq \epsilon_{\mathrm{p}}, \quad |C \bullet (U_t U_t^\top) - (-b^\top p_t - \tau^2\theta_t)|\leq \epsilon_{\mathrm{pd}}
\end{equation}
then \textbf{return} $(\bar U, \bar p, \bar\theta) = (U_t, p_t, \theta_{t})$;
\State set $t = t+1$ and \textbf{go to} step \textbf{3.}
\end{algorithmic}
\end{algorithm}
Several remarks about each of the steps of \ourmethod are now given. First, $t$ is the iteration count for \ourmethod. Second, in practice we choose the sequence $\{\beta_{t}\}_{t\geq 1}$ in an adaptive manner based on
some of the ideas of the LANCELOT method \cite{conn1991globally}.
Third, the while loop in lines 5 to 12 of \ourmethod corresponds to the HLR method, which is used by \ourmethod to approximately solve its AL subproblems in \eqref{eq:AL}. 
In line 4, the triple $(Y, G, \theta)$ is set as $(0_{n\times s},0_{n\times n},\infty)$ to ensure that the HLR method is always executed during each iteration of \ourmethod. Fourth, the nonconvex solver in line 6 of \ourmethod will typically make several evaluations of $\nabla_{U} \mathcal L_{\beta_t}(UU^{\top};p_{t-1})$, which is easily seen to be equal to
\begin{equation}\label{Grad U Space}
\nabla_{U} \mathcal L_{\beta_t}(UU^{\top};p_{t-1}):=2\left[C+\mathcal A^{*}\left(p_{t-1} +\beta_t\left(\mathcal A(UU^{\top})-b\right)\right)\right]U.
\end{equation}
Fifth, it follows from the definition of $\mathcal L_{\beta_t}(\cdot;p_{t-1})$ that $G$ in \eqref{gradient Def Algo 2} can be written as
\begin{equation}\label{Grad X Space}
    G:=C+\mathcal A^{*}\left(p_{t-1} +\beta_t\left(\mathcal A(YY^{\top})-b\right)\right).
\end{equation}
Sixth, after the HLR method approximately solves an AL subproblem, \ourmethod updates the Lagrange multiplier in line 14 and checks for termination in line 15. Lastly, the termination condition on this line does not check for near dual feasibility since it can be shown that for every $t \ge 1$,
\begin{align}\label{eq:St}
S_t := C+\mathcal A^*p_t + \theta_tI\succeq 0, \quad \theta_t\geq 0,
\end{align}
and hence that \ourmethod always generates feasible dual iterates
$(p_t,\theta_t,S_t)$.

We now comment on the steps performed within the while loop (lines 5 to 12) in light of the three steps of the \ourmethod outline given just before its formal description. Step 6 implements step i) of the outline. Steps 7 through 9 implement step ii). Finally, step 11 implements step iii).

We now give an interpretation of step 11 of \ourmethod relative to the $X$-space. The Frank-Wolfe procedure generates a new iterate by moving from the current point $YY^\top$ towards an extreme point of the feasible set $\Delta^n_{\tau}$, which is of the form $\tau^2 y y^\top$. The stepsize $\alpha$ computed in \eqref{FW stepsize} corresponds to solving
\[
\argmin_{\hat{\alpha} \in [0,1]}\left\{\mathcal L_{\beta_t}\left((1-\hat{\alpha})YY^{\top}+\hat{\alpha}(\tau yy^{\top});p_{t-1}\right) \right\}.
\]
Moreover, the iterate
$\tilde Y$ computed in 
equation \eqref{tilde Y Algo 2} has the property that
$\tilde Y \tilde Y^\top$ is equal to the usual Frank-Wolfe iterate $(1-\alpha)YY^{\top} + \alpha(\tau^2 yy^{\top})$ in the $X$-space; hence $\tilde Y$ is a valid $U$-factor for the Frank-Wolfe iterate.

\section{Computational Aspects of cuHALLaR}

GPU architecture is designed to perform many operations through massive parallelism, using a single instruction multiple data (SIMD) model. Unlike CPUs, which have fewer cores and are optimized for sequential execution and control, GPUs have thousands of simpler cores that are organized in 
streaming multiprocessors (SMs).
GPUs possess a memory hierarchy similar to CPUs with fast and small registers for each SM, a shared memory between blocks of SMs, and a slower and bigger global GPU memory.

In this section, we discuss the key computational factors that allow \gpuourmethod to achieve significant speedups over \ourmethod through the use of GPU programming.
This section is divided into two subsections. The first subsection discusses effectively dividing tasks between the CPU and the GPU. 
The second subsection describes our new way of efficiently reading and inputting SDP data that is suitable for GPU programming.

\subsection{Effective Division of CPU and GPU Tasks}
The primary computational concern in GPU programming is effectively dividing tasks between the CPU and GPU to exploit the advantages of both architectures.
Highly parallelizable operations, such as matrix, vector and tensor operations should be implemented in the GPU as functions which are called kernels. 
On the other hand, serial tasks should be implemented as regular functions in the CPU. 
Hence, the primary computational bottleneck in modern GPUs is the possibly large number of data transfers that need to be performed between the CPU and the GPU.

Another important computational consideration in GPU programming is that, like in the CPU, launching threads in the GPU also incurs an overhead due to thread management and communication.
Therefore, an efficient GPU implementation must minimize the number of data transfers and the number of kernels launched.

To minimize CPU and GPU communication overhead, \gpuourmethod pre-allocates essential problem data from \eqref{eq:sdp-primal} on the GPU, including the objective matrix $C$, the linear operator $\mathcal{A}$, and the right-hand side vector $b$. All subsequent computations involving large matrices and vectors are performed exclusively on the GPU.

\subsection{Reading Data Efficiently and Input Format}\label{Input format}

\paragraph{SDPA format.} Introduced in \cite{sdpa1,sdpa2,fujisawa2000sdpa}, the SDPA format allows for the specification of SDPs using sparse matrix representations of the cost and constraint matrices. Each matrix is stored by explicitly listing all nonzero entries in its upper triangular part, which is efficient for problems where all data matrices are sparse. However, for SDP instances with a dense but well-structured data matrix (e.g., the matrix of all ones or the identity matrix plus a low-rank matrix), it becomes prohibitively expensive to store and handle their input data in SDPA format. 
For example, if $C$ is the matrix of all ones, then SDPA format requires the storage of  $n(n+1)/2$ entries, and hence does not take advantage of the fact that such a matrix can be represented by either a single scalar or a single $n$-vector.

\paragraph{HSLR format.} To address the issues of the SDPA format and to improve handling of large-scale problems, we introduce a new input format, namely, Hybrid Sparse Low-Rank (HSLR) format, where each constraint matrix $A_\ell$ (with the convention $C = A_0$) can be decomposed as:
\begin{equation}\label{eq:HSLR-rep-body}
A_\ell = A_\ell^{\mathrm{sp}} + A_\ell^{\mathrm{lr}},
\qquad 
A_\ell^{\mathrm{lr}} = K_\ell D_\ell K_\ell^\top,
\qquad 
\ell=0,1,\dots,m
\end{equation}
where $A_\ell^{\mathrm{sp}}\in\mathbb{S}^n$ is a sparse symmetric matrix, $K_\ell\in\mathbb{R}^{n\times r_\ell}$, and $D_\ell\in\mathbb{S}^{r_\ell}$ is a symmetric matrix that is not necessarily diagonal. This structure allows us to reformulate the core operations of \gpuourmethod in terms of sparse matrix algebra and small dense matrix multiplications. In particular, \gpuourmethod never has to form any $n \times n$ matrices. For examples on how to construct HSLR format for several structured SDP instances such as the SDP relaxations of the matrix completion and maximum stable set problems, the reader should refer to our user manual~\cite{aguirre2025user}. 

The advantages of HSLR over sparse SDPA format are twofold:
\begin{enumerate}
    \item \textbf{Storage Efficiency}: 
    HSLR format has a significantly lower memory footprint than SDPA format. Specifically, if a data matrix $A_{\ell}$ is such that $A_{\ell}^{\mathrm lr}\neq 0$ and $D_{\ell}$ is diagonal, then it requires $n_{\ell}+r_{\ell}n+r_{\ell}$ entries (resp., n(n+1)/2) to be stored in HSLR (resp., SDPA) format, where $n_{\ell}$ is the number of nonzeros of $A_{\ell}^{\mathrm sp}$.
    
    \item \textbf{Computational Efficiency}: During \gpuourmethod's execution, matrix-vector products and trace evaluations are computed by separately handling the sparse and low-rank components of each data matrix $A_{\ell}$. Since \gpuourmethod handles these components separately, it never has to form the possibly dense matrix $A_{\ell}$ and calls different kernels based on the specific structure of each $A_\ell$.
\end{enumerate}

\section{Efficient Implementation of  $\tilde {\mathcal A}(UU^{\top})$ and $\tilde {\mathcal A}^{*}(\tilde{p})U$ }\label{sec:gpu_ops}

The efficient evaluations of the linear maps $\mathcal A(X):\mathbb S^{n}\rightarrow \mathbb R^{m}$ and $\mathcal A^{*}(p):\mathbb R^{m}\rightarrow \mathbb S^{n}$ as in \eqref{Definition of Linear Operators} are essential to the performance of \gpuourmethod. To avoid the formation of a possibly dense $n\times n$ matrix $X$, \gpuourmethod forms a low-rank factor $U\in\mathbb R^{n\times r}$ of $X$ that satisfies $X=UU^{\top}$. At every iteration, \gpuourmethod has to evaluate for a given $(U, p) \in \mathbb R^{n\times r} \times \R^m$ the quantities  $\tilde{\cA}(UU^\top)$ and $\tilde{\cA}^*(\tilde{p})U$ where $\tilde p=(1,p) \in \R^{m+1}$, and
$\tilde {\mathcal A}:\mathbb S^{n}\rightarrow \mathbb R^{m+1}$ and $\tilde{\cA}^*:\mathbb{R}^{m+1}\to\mathbb{S}^n$ are given by
\begin{equation}\label{A tilde operator}
\tilde {\mathcal A}(X):=\begin{bmatrix}
           A_0\bullet X \\
           A_1\bullet X \\
           \vdots \\
           A_m \bullet X
         \end{bmatrix}:=\begin{bmatrix}
           A_0\bullet X \\
           \mathcal A(X)
         \end{bmatrix}\,,\qquad \tilde {\mathcal A}^{*}(\tilde{p}):=\sum_{\ell=0}^{m}\tilde p_\ell A_{\ell}=A_0+\sum_{\ell=1}^{m}p_\ell A_{\ell} = A_0 + 
         {\mathcal A}^{*}( p)
\end{equation}
where $A_{0}:=C$ and $\mathcal A(\cdot)$ and $\mathcal A^{*}(\cdot)$ are as in \eqref{Definition of Linear Operators}. It follows from the definitions of  $\nabla_{U} \mathcal L_{\beta_t}(UU^{\top};p_{t-1})$ and $G$ in \eqref{Grad U Space} and \eqref{Grad X Space}, respectively, that \gpuourmethod has to evaluate $\tilde {\mathcal{A}}(UU^{\top})$ in steps 6, 7, 11, and 14, and $\tilde {\mathcal A}^{*}(\tilde{p})U$ in steps 6 and 11.

This section presents the details of the efficient implementation of the operations in \eqref{A tilde operator} in the GPU, using the fact that the data matrices $A_{\ell}$ are assumed to be the sum of a sparse and a low-rank matrix
as in~\eqref{eq:HSLR-rep-body}.

\subsection{Constraint Evaluation $\tilde {\mathcal A}(UU^{\top})$}\label{App:ConstraintEval}

This subsection describes how \gpuourmethod efficiently computes $\tilde {\mathcal A}(UU^{\top})$, where $U\in \mathbb R^{n\times r}$ and $\tilde A(\cdot)$ is as in \eqref{A tilde operator}. It follows from \eqref{eq:HSLR-rep-body} and the first relation in \eqref{Definition of Linear Operators} that $T=\tilde {\mathcal A}(UU^{\top})$ can be computed as $T=T^{\mathrm sp}+T^{\mathrm lr}$ where the $\ell$-th components of $T^{\mathrm sp}$ and $T^{\mathrm lr}$ are defined as:
\begin{equation}
T^{\mathrm sp}_\ell =A^{\mathrm sp}_\ell \bullet (UU^\top), \quad T_{\ell}^{\mathrm lr}=A^{\mathrm lr}_\ell \bullet (UU^\top), \quad \ell=0,1,\dots,m.
\end{equation}
In practice, \gpuourmethod implements efficient subroutines to separately compute $T^{\mathrm sp}$ and $T^{\mathrm lr}$.
The rest of this subsection discusses in more detail how \gpuourmethod computes these two quantities efficiently.

\paragraph{Sparse Component.}
This paragraph describes how \gpuourmethod efficiently computes $T^{\mathrm sp}\in \mathbb R^{m+1}$. 
Let $B=UU^{\top}$. It is easy to see that the entries of $B$ can be computed as $B_{ij}=\inner{U_{i,:}}{U_{j,:}}$, where $U_{i,:}$ and $U_{j,:}$ are the $i$-th and $j$-th rows of $U$, respectively.
The $\ell$-th component of $T^{\mathrm{sp}}$, $T_\ell^{\mathrm{sp}}=A_\ell^{\mathrm{sp}} \bullet (UU^\top)$, is then computed as:
\begin{equation}
T_\ell^{\mathrm{sp}} = \sum_{i=1}^n \sum_{j=1}^n (A_\ell^{\mathrm{sp}})_{ij} B_{ij}.
\end{equation}
Since $A_\ell^{\mathrm{sp}}$ and $B$ are symmetric matrices, this summation can be rewritten by defining the weights
\begin{equation}\label{omega}
\omega_{ij,\ell} \defi (2-\delta_{ij})(A_\ell^{\mathrm{sp}})_{ij}, \text{ where } \delta_{ij}=
\begin{cases}
1 \quad \text{if } i=j\\
0 \quad \text{otherwise}.
\end{cases}
\end{equation}
Using \eqref{omega}, we see that $T_\ell^{\mathrm{sp}} = \sum_{i \le j} \omega_{ij,\ell} B_{ij}$.

To avoid redundant computations of $B_{ij}$ across different constraints $\ell$, we reorganize this computation by using the union of the supports of $A_{\ell}^{\mathrm sp}$. This union $\mathcal S$ is formally defined as:
\begin{equation}\label{def: support S ell}
\mathcal{S}\defi\bigcup_{\ell=0}^m \mathrm{supp}(A_\ell^{\mathrm{sp}})\cap \{(i,j):1\le i\le j\le n\}.
\end{equation}
For each pair $(i,j) \in \mathcal{S}$, we define $\mathcal L_{ij}$ to be the set of indices $\ell$ for which the entry $(i,j)$ is active:
\begin{equation}\label{eq: index list}
\mathcal{L}_{ij}\defi\{\ell:(i,j)\in \mathrm{supp}(A_\ell^{\mathrm{sp}})\}.
\end{equation}
To compute $T^{\mathrm sp}$, we then compute the quantity $B_{ij}$ over each pair $(i,j)\in \mathcal S$ and update the corresponding entries $T^{\mathrm sp}_\ell$ for all $\ell \in \mathcal{L}_{ij}$.
This procedure is formally described in Algorithm~\ref{alg:A_UUT_sparse} below.

\begin{algorithm}[H]
\setstretch{0.8}
\caption{\normalsize Computation of $\tilde {\mathcal A}(UU^{\top})$ for the sparse component, where $\tilde {\mathcal A}(\cdot)$ is as in \eqref{A tilde operator}.}
\label{alg:A_UUT_sparse}
\begin{algorithmic}[1]
\normalsize
\Require $U \in \mathbb{R}^{n \times r}$, data matrices $\{A_\ell^{\mathrm{sp}}\}_{\ell=0}^m \subseteq \mathbb{S}^n$,  
and precomputed structures $\mathcal{S}$ and $\{\mathcal{L}_{ij}\}$ from~\eqref{def: support S ell}–\eqref{eq: index list}.
\Ensure $T^{\mathrm sp} \in \mathbb{R}^{m+1}$ where $T^{\mathrm sp}_\ell = A^{\mathrm sp}_{\ell} \bullet UU^\top$.
\vspace{3pt}
\State Initialize $T^{\mathrm sp} \gets 0 \in \mathbb{R}^{m+1}$.
\vspace{3pt}
\ForAll{$(i,j) \in \mathcal{S}$ in parallel}
    \State $B_{ij} \gets \inner{U_{i,:}}{U_{j,:}}$
    \ForAll{$\ell \in \mathcal{L}_{ij}$}
        \State 
        $\omega_{ij,\ell} \gets (2-\delta_{ij})(A_\ell^{\mathrm{sp}})_{ij}$
        \vspace{4pt}
        \State $T^{\mathrm sp}_{\ell} \gets T^{\mathrm sp}_{\ell} + \omega_{ij,\ell} B_{ij}$
        \vspace{4pt}
    \EndFor
\EndFor
\Statex
\vspace{2pt}
\Return{$T^{\mathrm sp}$}
\end{algorithmic}
\end{algorithm}

In Algorithm~\ref{alg:A_UUT_sparse}, the outer loop over the set $\mathcal{S}$ is parallelized. More specifically, each unique index pair $(i,j) \in \mathcal{S}$ is assigned to a separate GPU thread. Each thread performs a loop over $r \ll n$ to compute $B_{ij}$ and a loop over $|\mathcal{L}_{ij}| \ll m+1$ that updates $T^{\mathrm sp}_\ell$.

\paragraph{Low-Rank Component.}
This paragraph describes how \gpuourmethod efficiently computes $T^{\mathrm lr}\in \mathbb R^{m+1}$. 
Using the second relation in \eqref{eq:HSLR-rep-body}, the $\ell$-th component of $T^{\mathrm{lr}}$,
$T_\ell^{\mathrm{lr}}=A_\ell^{\mathrm{lr}} \bullet (UU^\top)$,
can be computed as
\begin{equation}\label{T LR}
T_\ell^{\mathrm{lr}} 
=\mathrm{Tr}(K_\ell D_\ell K_\ell^\top (U U^\top)) = \mathrm{Tr}(D_\ell (K_\ell^\top U) (U^\top K_\ell)).
\end{equation}
Defining $Y_\ell \defi K_\ell^\top U \in \mathbb{R}^{r_\ell\times r}$, relation \eqref{T LR} further simplifies to:
\begin{equation}\label{eq:map-lr-appA}
T_\ell^{\mathrm{lr}} = \mathrm{Tr}(D_\ell Y_\ell Y_\ell^\top).
\end{equation}
Algorithm~\ref{alg:A_UUT_low_rank} below formally describes our efficient procedure for computing $T^{\mathrm{lr}}$.
\begin{algorithm}[H]
\setstretch{0.8}
\caption{\normalsize Computation of $\tilde {\mathcal A}(UU^{\top})$ for the low-rank component, where $\tilde {\mathcal A}(\cdot)$ is as in \eqref{A tilde operator}}
\label{alg:A_UUT_low_rank}
\begin{algorithmic}[1]
\normalsize
\Require $U \in \mathbb{R}^{n \times r}$ and problem data $\{K_\ell, D_\ell\}_{\ell=0}^m$.
\Ensure $T^{\mathrm lr} \in \mathbb{R}^{m+1}$ where $T^{\mathrm lr}_\ell = A^{\mathrm lr}_{\ell} \bullet UU^\top$.
\State Initialize $T^{\mathrm lr} \gets 0 \in \mathbb{R}^{m+1}$.
\For{$\ell=0$ to $m$}
    \vspace{2pt}
    \State $Y_\ell \gets K_\ell^\top U$ \Comment{in parallel}
    \vspace{2pt}
    \State $T^{\mathrm lr}_\ell \gets 0 + \mathrm{Tr}(D_\ell Y_\ell Y_\ell^\top)$
\EndFor
\Statex
\Return{$T^{\mathrm lr}$}
\end{algorithmic}
\end{algorithm}
In Algorithm~\ref{alg:A_UUT_low_rank}, the loop over $\ell=0,\dots,m$ is performed in serial, while the computations of $Y_\ell$ and $\mathrm{Tr}(D_\ell Y_\ell Y_\ell^\top)$ are performed in parallel. In our numerical experiments, we found that parallelizing these computations was very efficient.

\subsection{Adjoint $\tilde {\mathcal A}^{*}(\tilde{p})U$}\label{App:Adjoint}

This subsection details how, for a given $(U, p) \in \mathbb R^{n\times r} \times \R^m$, \gpuourmethod efficiently computes 
the quantity $W=\tilde {\mathcal A}^{*}(\tilde{p})U$ as in \eqref{A tilde operator}, where $\tilde p=(1,p)\in \mathbb R^{m+1}$. It follows from \eqref{eq:HSLR-rep-body} and the second identity in \eqref{A tilde operator} that $W$ is given by
\begin{equation}\label{Adjoint calculation}
W = \underbrace{\left(\sum_{\ell=0}^m \tilde p_\ell A_\ell^{\mathrm{sp}}\right)U}_{W^{\mathrm{sp}}} + \underbrace{\sum_{\ell=0}^m \tilde p_\ell (K_\ell D_\ell K_\ell^\top) U}_{W^{\mathrm{lr}}}\,.
\end{equation}
The remaining part of this subsection discusses how \gpuourmethod implements efficient subroutines to compute $W^{\mathrm sp}$ and $W^{\mathrm lr}$ separately.

\paragraph{Sparse Component.}
To compute $W^{\mathrm sp}$ as in \eqref{Adjoint calculation}, we first allocate a sparse symmetric matrix $M$ with nonzero entries at indices $(i,j)$, where $\mathcal{L}_{ij}$ in \eqref{eq: index list} is nonempty. After doing this allocation, we then update the values of the nonzero entries of $M_{ij}$ in parallel using the formula:
\begin{equation}
M_{ij} = \sum\limits_{\ell \in \mathcal{L}_{ij}} \tilde p_\ell (A^{\mathrm{sp}}_{\ell})_{ij}.
\end{equation}
The matrix $W^{\mathrm{sp}}$ is then computed as $W^{\mathrm{sp}}=MU$. This procedure is formally described in Algorithm~\ref{alg:Astar_lambda_U_sp} below.

\begin{algorithm}[H]
\setstretch{0.8}
\caption{\normalsize Computation of $W^{\mathrm sp}$ as in \eqref{Adjoint calculation}.}
\label{alg:Astar_lambda_U_sp}
\begin{algorithmic}[1]
\normalsize
\Require $U \in \mathbb{R}^{n \times r}$, $\tilde p=(1,p) \in \mathbb{R}^{m+1}$ where $p\in \mathbb R^{m}$, $\{A_\ell^{\mathrm{sp}}\}_{\ell=0}^m \subseteq \mathbb{S}^n$,  
and precomputed structures $\mathcal{S}$, $\{\mathcal{L}_{ij}\}$ from~\eqref{def: support S ell}–\eqref{eq: index list}.
\Ensure $W^{\mathrm sp} \in \mathbb{R}^{n \times r}$ where $W^{\mathrm sp} = \left(\sum_{\ell=0}^m \tilde p_\ell A_\ell^{\mathrm{sp}}\right)U$.
\State Initialize $M \in \mathbb{S}^{n}$ with pattern $\mathcal S\cup\{(j,i):(i,j)\in\mathcal S\}$.
\ForAll{$(i,j) \in \mathcal{S}$, in parallel}
    \State $w \gets 0$
    \ForAll{$\ell \in \mathcal{L}_{ij}$}
        \State $w \gets w + \tilde p_{\ell}(A^{\mathrm{sp}}_\ell)_{ij}$
    \EndFor
    \State $M_{ij} \gets w$
    \If{$i \neq j$}
        \State $M_{ji} \gets w$
    \EndIf
\EndFor
\State $W^{\mathrm sp} \gets MU$
\Statex
\Return{$W^{\mathrm sp}$}
\end{algorithmic}
\end{algorithm}
Several remarks about Algorithm~\ref{alg:Astar_lambda_U_sp} are now given. The matrix $M$, which is initialized in step 3, is stored in sparse format with symmetric pattern $\mathcal S\cup\{(j,i):(i,j)\in\mathcal S\}$, hence it contains at most $2|\mathcal S| - |\{\,i : (i,i)\in\mathcal S\,\}|$ nonzero entries. Also, it is worthwhile to note that the outer for loop over $\mathcal{S}$ is parallelized. More specifically, each unique index pair $(i,j) \in \mathcal{S}$ is assigned to a separate GPU thread. Each thread then performs a loop over $|\mathcal{L}_{ij}| \ll m+1$ to compute  $w = \sum_{\ell \in \mathcal{L}_{ij}} \tilde p_{\ell}(A^{\mathrm{sp}}_\ell)_{ij}$ and updates the corresponding entries in $M$.

\paragraph{Low-Rank Component.}
This paragraph describes how \gpuourmethod efficiently computes $W^{\mathrm lr}\in \mathbb R^{n \times r}$. Using the second quantity in \eqref{Adjoint calculation} and defining
$Y_\ell = K_\ell^\top U \in \mathbb{R}^{r_\ell\times r}$, it follows that $W^{\mathrm{lr}}$ can be computed as
\begin{equation}
W^{\mathrm{lr}} = \sum_{\ell=0}^m \tilde p_\ell K_\ell (D_\ell Y_\ell).
\end{equation}
Algorithm~\ref{alg:Astar_lambda_U_lr} below formally describes our efficient procedure for computing $W^{\mathrm{lr}}$.

\begin{algorithm}[H]
\setstretch{0.8}
\caption{\normalsize Computation of $W^{\mathrm lr}$ as in \eqref{Adjoint calculation}.}
\label{alg:Astar_lambda_U_lr}
\begin{algorithmic}[1]
\normalsize
\Require $U \in \mathbb{R}^{n \times r}$, $\tilde p=(1,p) \in \mathbb{R}^{m+1}$ where $p\in \mathbb R^{m}$,
and problem data $\{K_\ell, D_\ell\}_{\ell=0}^{m}$.
\Ensure $W^{\mathrm lr} \in \mathbb{R}^{n \times r}$ where $W^{\mathrm lr} = \left(\sum_{\ell=0}^m \tilde p_\ell A_\ell^{\mathrm{lr}}\right)U$.
\State Initialize $W^{\mathrm lr} \gets 0_{n\times r}$.
\vspace{3pt}
\For{$\ell=0$ to $m$}
\vspace{3pt}
    \State $Y_\ell \gets K_\ell^\top U$
    \vspace{3pt}
    \State $Z_\ell \gets D_\ell Y_\ell$ \Comment{in parallel}
    \vspace{3pt}
    \State $W^{\textrm lr} \gets W^{\textrm lr} + \tilde p_\ell (K_\ell Z_\ell)$
    \vspace{3pt}
\EndFor
\Statex
\Return{$W^{\textrm lr}$}
\end{algorithmic}
\end{algorithm}
In Algorithm~\ref{alg:Astar_lambda_U_lr}, the loop over $\ell=0,\dots,m$ is performed in serial, while the computations of $Y_\ell = K_\ell^\top U$ and $Z_\ell = D_\ell Y_\ell$ used to update $W^{\mathrm lr}$ are performed in parallel.

\subsection{Kernels for problems in SDPA format}\label{subsec: SDPA kernel details}
This subsection describes \gpuourmethod's efficient approach for dealing with SDP instances given in sparse SDPA format. 
The parser reads each data matrix $A_{\ell}$, for $\ell=0,...,m$ and computes its density as
\begin{equation}\label{eq: constraint matrix density}
d_\ell = \frac{nnz(A_\ell)}{n^2}.
\end{equation}
The kernels for SDPA format separate constraints into sparse or dense. Let $A_\ell^{dn}$ be the representation of data matrix $A_\ell$ in dense format. If $d_\ell > 10^{-1}$, \gpuourmethod sets $A_\ell = A_\ell^{dn}$; otherwise it sets $A_\ell = A_\ell^{\mathrm sp}$. The quantities $\mathcal{S}$ and $\{\mathcal{L}_{ij}\}$ in \eqref{def: support S ell} and \eqref{eq: index list} are precomputed from the data matrices which are determined to be sparse. The evaluations of $\tilde {\mathcal A}(UU^{\top})$ and $\tilde {\mathcal A}^{*}(\tilde{p})U$ use similar approaches as the ones described in Section~\ref{sec:gpu_ops}. The sparse components of these evaluations are performed using Algorithms~\ref{alg:A_UUT_sparse} and~\ref{alg:Astar_lambda_U_sp} while the dense components of these evaluations are evaluated using adaptations of Algorithms~\ref{alg:A_UUT_low_rank} and~\ref{alg:Astar_lambda_U_lr}, where  $A_\ell^{dn}$ replaces $A_\ell^{\mathrm lr} = K_\ell D_\ell K_\ell^{\top}$.

\section{Numerical Experiments}

In this section, we provide extensive experiments comparing \gpuourmethod, \ourmethod, and \culorads on three problem classes: matrix completion, maximum stable set, and phase retrieval. These problems were also used in the benchmark of the original \ourmethod paper \cite{monteiro2024low}. The pre-compiled binary solver can be found at \url{https://github.com/OPTHALLaR}.

We attempted to run CuClarabel \cite{chen2024cuclarabel} on the matrix completion and stable set experiments, but CuClarabel failed even on relatively small SDPA files. Moreover, CuClarabel's lack of native SDPA support precluded straightforward benchmarking, so we omit comparisons against it.

The comparisons between \gpuourmethod and \culorads evaluate \gpuourmethod using the HSLR format described in Section~\ref{sec:gpu_ops} against \culorads using SDPA format. This comparison highlights the performance and memory efficiency advantages of our HSLR representation and \gpuourmethod's implementation.

\subsection{Experimental Setup}
\paragraph{Hardware and software. }
Our experiments compare the GPU methods \gpuourmethod and \culorads,
and the CPU method \ourmethod.
The GPU experiments are performed on an NVIDIA H200 with 142 GB VRAM deployed on a cluster with an Intel Xeon Platinum 8469C CPU. For CPU benchmarks, we use a single Dual Intel Xeon Gold 6226 CPUs @ 2.7 GHz (24 cores/node).
Further details are given in Table~\ref{tab:flops}. We set the memory of the CPU benchmarks at 142GB, which is equal to the maximum memory of the GPU.

We implemented \gpuourmethod as a Julia \cite{bezanson2017julia} module. Our implementation uses the CUDA platform \cite{choquette2023nvidia} for interfacing with NVIDIA CUDA GPUs. The experiments are performed using CUDA 12.6.1 and Julia 1.11.3. The reported running times do not take into account the pre-compilation time.

\begin{table}[ht!]
\centering
\small
\begin{tabular}{rll}
\toprule
\textbf{Specification}      & \textbf{CPU Node}                                        & \textbf{GPU}                               \\
\midrule
\textbf{Processor(s)}       & Dual Intel Xeon Gold 6226 @ 2.70\,GHz   & NVIDIA H200                                \\
\textbf{Cores / SMs}        & 24 cores                           & 80 SMs                                    \\
\textbf{Cache (per CPU)}    & L3: 19.25\,MB                                            &                                           \\
\textbf{\ (per core)}       & L2: 1\,MB; L1: 32\,KB (D+I)                              & L1/Shared (per SM), L2             \\
\textbf{Peak FP64 Perf.}    & $\sim$2.07\,TFLOPS                        & $\sim$34\,TFLOPS                  \\
\textbf{Memory Size}        & 142\,GB DDR4                                             & 142\,GB HBM3e                              \\
\textbf{Memory Bandwidth}   & $\sim$140.7\,GB/s                & 4.8\,TB/s                                  \\
\bottomrule
\end{tabular}
\caption{Comparison of CPU and GPU specifications used in experiments.}
\label{tab:flops}
\end{table}

\FloatBarrier

\paragraph{Termination criteria.}
Given a tolerance $\epsilon > 0$, both solvers stop when a set of relative error metrics falls below $\epsilon$. For \gpuourmethod, an iterate $(X, p, \theta) \in \Delta^n_{\tau} \times \mathbb R^{m} \times \mathbb R_+$ is deemed optimal if it satisfies
\begin{align}\label{eq:termination2}
\max \left\{
\frac{\|\mathcal{A}X-b\|_2}{1+\|b\|_1},
\frac{|\mathrm{pval} - \mathrm{dval}|}{1+|\mathrm{pval}|+|\mathrm{dval}|},
\frac{\max\{0,-\lambda_{\min}(S) \}}{1+\|C\|_1}
\right\} \leq \epsilon,
\end{align}
where $\mathrm{pval} = C \bullet X$ is the primal value, $S = C+{\mathcal{A}}^{*}{p} + \theta I$ is the dual slack matrix, and the dual value is $\mathrm{dval}=-b^\top p - \tau^2\theta$. The solver \culorads employs a similar set of metrics; however, as it addresses a formulation without the spectraplex constraint, its dual value is defined as $\mathrm{dval} = -b^\top p$.

\paragraph{Initialization.}
Both \gpuourmethod and \ourmethod are initialized by setting the dual variable to $p_0 = 0$ and the primal factor $U_0$ to an $n \times 1$ matrix  with entries drawn independently from a standard Gaussian distribution.

\paragraph{Input format.} \gpuourmethod binary accepts two different input formats: HSLR and SDPA formats. For SDP data prepared in SDPA format, the user must also supply the trace bound $\tau > 0$ in \eqref{Delta Definition} for the constraint $\mathrm{Tr}(X) \leq \tau^2$ externally via the command-line interface or a parameter file. On the other hand, HSLR embeds the trace bound $\tau$ directly within the file itself, so that the user does not need to supply it separately. For more details and examples of preparing data in SDPA and HSLR format, the reader should refer to our user manual \cite{aguirre2025user}.

Our computational experiments also report results for \ourmethod and \gpuourmethod versions (not available to the public) where the data is generated in the memory. Specifically, instead of reading the problem data from an input file, these versions use subroutines for computing $\tilde{\mathcal{A}}(UU^\top)$ and $\tilde{\mathcal{A}}^*(\tilde p)U$ for any $\tilde p=(1,p)\in\mathbb{R}^{m+1}$ and $U\in\mathbb{R}^{n\times r}$. Moreover, these two ``on-the-fly'' versions have the advantage of being able to run the phase retrieval experiments, which are expensive to store in either SDPA or HSLR.

\subsection{Experiments for matrix completion}
Given integers $n_2 \geq n_1 \ge 1$, consider the problem of retrieving a low rank matrix $M \in \mathbb{R}^{n_1 \times n_2}$ by observing a subset $\{ M_{ij}:  (i,j) \in \Omega\}$ of its entries. A standard approach to tackle this problem is by considering the nuclear norm relaxation:
\[
\min_{Y \in \mathbb{R}^{n_1\times n_2}} \quad \left\{\|Y\|_*
\ : \ Y_{ij} = M_{ij}, \ \forall\, (i,j) \in \Omega
\right \}.
\]
The above problem can be rephrased as the following SDP:
\begin{align}
\label{eq:matcomp-sdp}
\min_{X \in \mathbb S^{n_1+n_2}}
\left\{ \,\frac12 \mathrm{Tr}(X) :
X= \begin{pmatrix} W_1 & Y \\ Y^\top & W_2 \end{pmatrix} \succeq 0\,, \quad
Y_{i,j} = M_{i,j}\,, \; \forall\, (i,j) \in \Omega)
\right\}.
\end{align}
Matrix completion instances are generated randomly, using the following procedure: given $r\le n_1\le n_2$, the hidden solution matrix $M$ is the product $UV^\top$, where the matrices $U\in\mathbb{R}^{n_1\times r}$ and $V\in\mathbb{R}^{n_2\times r}$ have independent standard Gaussian random variables as entries. Afterward, $m$ independent and uniformly random entries from $M$ are taken, where $m=\lceil \gamma \ r(n_1+n_2-r)\rceil$, and $\gamma=r\log(n_1+n_2)$ is the oversampling ratio.

Table~\ref{tab:gpuourmethod_vs_cpu} compares \ourmethod, \gpuourmethod, and \culorads on matrix completion problems. The second block of columns demonstrates the tremendous speedups our GPU-accelerated solver \gpuourmethod achieves over \ourmethod using on-the-fly-problem generation. On many instances, \gpuourmethod achieves speedups of 10x to 188x compared to \ourmethod.
For the largest instance considered with $(n_1, n_2) = (800,\!000, 1,\!200,\!000)$, \gpuourmethod solves the problem in 53 seconds, while \ourmethod requires nearly 2.5 hours.

The last block of columns compares \gpuourmethod, using both HSLR and SDPA formats via our binary, against \culorads, which uses SDPA format. Here, \gpuourmethod with HSLR format consistently outperforms \culorads across all instances, achieving speedups ranging from 0.9x to 37x. Notably, \gpuourmethod efficiently finds near-optimal solutions with much smaller ranks than the ones found \culorads.
For the largest instances considered in the table below the midrule line, \gpuourmethod with HSLR format achieves 2x–4x speedups over \culorads. 

\FloatBarrier

\begin{table}[!tbh]
\captionsetup{font=scriptsize}
\begin{centering}
\begin{tabular}{>{\centering\arraybackslash}p{4.2cm} >{\centering\arraybackslash}p{0.9cm} | >{\centering\arraybackslash}p{1.5cm} >{\centering\arraybackslash}p{1.5cm} >{\centering\arraybackslash}p{1.2cm} || >{\centering\arraybackslash}p{1.1cm} >{\centering\arraybackslash}p{1.2cm} >{\centering\arraybackslash}p{1.2cm} >{\centering\arraybackslash}p{0.75cm}}
\multicolumn{2}{c}{\textbf{\scriptsize{Problem Instance}}} & \multicolumn{3}{c}{\textbf{\scriptsize{\ourmethod vs \gpuourmethod (on-the-fly)}}} & \multicolumn{4}{c}{\textbf{\scriptsize{\gpuourmethod vs \culorads}}} \tabularnewline
\toprule
\scriptsize{Size $(n_1, n_2 ; m)$} & \scriptsize{$r$} & \scriptsize{CPU} &\scriptsize{GPU} & \scriptsize{Ratio} & \scriptsize{\gpuourmethod (HSLR)} & \scriptsize{\gpuourmethod (SDPA)} & \scriptsize{\culorads (SDPA)} & \scriptsize{Ratio$^*$} \tabularnewline
\midrule
\scriptsize{3,000, 7,000; 828,931}    & \scriptsize{3} & \scriptsize{7.63/3} & \scriptsize{0.77/3} & \scriptsize{9.9} & \scriptsize{3.33/3} & \scriptsize{3.22/3} & \scriptsize{123.71/19} & \scriptsize{ {\textbf{\color{blue}37.1}}  } \tabularnewline
\scriptsize{3,000, 7,000; 828,931}    & \scriptsize{3} & \scriptsize{7.52/3} & \scriptsize{0.70/3} & \scriptsize{10.7} & \scriptsize{3.20/3} & \scriptsize{3.25/3} & \scriptsize{11.42/19} & \scriptsize{{\textbf{\color{blue}3.6}}} \tabularnewline
\scriptsize{3,000, 7,000; 2,302,586}    & \scriptsize{5} & \scriptsize{34.52/5} & \scriptsize{1.28/5} & \scriptsize{27.0} & \scriptsize{3.37/5} & \scriptsize{2.49/5} & \scriptsize{3.19/19} & \scriptsize{0.9} \tabularnewline
\scriptsize{3,000, 7,000; 2,302,586}    & \scriptsize{5} & \scriptsize{34.42/5} & \scriptsize{1.09/5} & \scriptsize{31.6} & \scriptsize{2.45/5} & \scriptsize{2.68/5} & \scriptsize{2.26/19} & \scriptsize{0.9} \tabularnewline
\scriptsize{10,000, 21,623; 2,948,996}  & \scriptsize{3} & \scriptsize{38.88/3} & \scriptsize{1.22/3} & \scriptsize{31.9} & \scriptsize{3.31/3} & \scriptsize{2.68/3} & \scriptsize{5.14/21} & \scriptsize{{\textbf{\color{blue}1.5}}} \tabularnewline
\scriptsize{10,000, 21,623; 2,948,996}  & \scriptsize{3} & \scriptsize{34.80/3} & \scriptsize{1.11/3} & \scriptsize{31.4} & \scriptsize{4.10/3} & \scriptsize{2.94/3} & \scriptsize{5.65/21} & \scriptsize{{\textbf{\color{blue}1.4}}} \tabularnewline
\scriptsize{10,000, 21,623; 8,191,654}  & \scriptsize{5} & \scriptsize{141.17/5} & \scriptsize{2.05/5} & \scriptsize{68.9} & \scriptsize{4.02/5} & \scriptsize{7.38/5} & \scriptsize{4.71/21} & \scriptsize{{\textbf{\color{blue}1.2}}} \tabularnewline
\scriptsize{10,000, 21,623; 8,191,654}  & \scriptsize{5} & \scriptsize{140.55/5} & \scriptsize{1.95/5} & \scriptsize{72.1} & \scriptsize{4.02/5} & \scriptsize{7.74/5} & \scriptsize{9.00/21} & \scriptsize{{\textbf{\color{blue}2.2}}} \tabularnewline
\scriptsize{25,000, 50,000; 3,367,574}  & \scriptsize{2} & \scriptsize{43.52/2} & \scriptsize{1.20/2} & \scriptsize{36.3} & \scriptsize{4.40/2} & \scriptsize{3.96/2} & \scriptsize{10.77/23} & \scriptsize{{\textbf{\color{blue}2.4}}} \tabularnewline
\scriptsize{25,000, 50,000; 7,577,040}  & \scriptsize{3} & \scriptsize{79.93/3} & \scriptsize{1.82/3} & \scriptsize{43.9} & \scriptsize{4.19/3} & \scriptsize{3.68/3} & \scriptsize{21.31/23} & \scriptsize{{\textbf{\color{blue}5.1}}} \tabularnewline
\scriptsize{30,000, 70,000; 4,605,171}  & \scriptsize{2} & \scriptsize{60.57/2} & \scriptsize{1.63/2} & \scriptsize{37.2} & \scriptsize{4.43/2} & \scriptsize{3.48/2} & \scriptsize{17.02/24} & \scriptsize{{\textbf{\color{blue}3.8}}} \tabularnewline
\scriptsize{30,000, 70,000; 10,361,633}  & \scriptsize{3} & \scriptsize{145.76/3} & \scriptsize{2.57/3} & \scriptsize{56.7} & \scriptsize{5.12/3} & \scriptsize{4.54/3} & \scriptsize{18.45/24} & \scriptsize{{\textbf{\color{blue}3.6}}} \tabularnewline
\scriptsize{50,000, 100,000; 7,151,035} & \scriptsize{2} & \scriptsize{121.63/2} & \scriptsize{1.89/2} & \scriptsize{64.4} & \scriptsize{12.36/2} & \scriptsize{5.18/2} & \scriptsize{23.67/24} & \scriptsize{{\textbf{\color{blue}1.9}}} \tabularnewline
\scriptsize{50,000, 100,000; 16,089,828} & \scriptsize{3} & \scriptsize{242.00/3} & \scriptsize{3.68/3} & \scriptsize{65.8} & \scriptsize{7.04/3} & \scriptsize{5.69/3} & \scriptsize{27.95/24} & \scriptsize{{\textbf{\color{blue}4.0}}} \tabularnewline
\scriptsize{80,000, 120,000; 9,764,859} & \scriptsize{2} & \scriptsize{237.38/2} & \scriptsize{2.45/2} & \scriptsize{96.9} & \scriptsize{5.63/2} & \scriptsize{6.05/2} & \scriptsize{20.66/25} & \scriptsize{{\textbf{\color{blue}3.7}}} \tabularnewline
\scriptsize{80,000, 120,000; 21,970,931} & \scriptsize{3} & \scriptsize{490.07/3} & \scriptsize{5.02/3} & \scriptsize{97.6} & \scriptsize{6.92/3} & \scriptsize{6.27/3} & \scriptsize{32.28/25} & \scriptsize{{\textbf{\color{blue}4.7}}} \tabularnewline
\scriptsize{120,000, 180,000; 34,051,152}& \scriptsize{3} & \scriptsize{1023.0/4} & \scriptsize{7.35/4} & \scriptsize{139.2} & \scriptsize{9.85/4} & \scriptsize{17.62/4} & \scriptsize{41.24/26} & \scriptsize{{\textbf{\color{blue}4.2}}} \tabularnewline
\scriptsize{160,000, 240,000; 46,437,192}& \scriptsize{3} & \scriptsize{1419.26/4} & \scriptsize{8.79/4} & \scriptsize{161.5} & \scriptsize{12.24/3} & \scriptsize{12.12/3} & \scriptsize{40.65/26} & \scriptsize{{\textbf{\color{blue}3.3}}} \tabularnewline
\scriptsize{240,000, 360,000; 71,845,299}& \scriptsize{3} & \scriptsize{1949.8/3} & \scriptsize{12.82/3} & \scriptsize{152.1} & \scriptsize{28.63/3} & \scriptsize{20.76/3} & \scriptsize{67.11/27} & \scriptsize{{\textbf{\color{blue}2.3}}} \tabularnewline
\midrule
\scriptsize{320,000, 480,000; 97,865,043}& \scriptsize{3} & \scriptsize{2754.83/3} & \scriptsize{17.41/3} & \scriptsize{158.2} & \scriptsize{36.76/3} & \scriptsize{31.78/3} & \scriptsize{136.43/28} & \scriptsize{{\textbf{\color{blue}3.7}}} \tabularnewline
\scriptsize{400,000, 600,000; 124,339,596}& \scriptsize{3} & \scriptsize{3472.13/3} & \scriptsize{22.31/3} & \scriptsize{155.6} & \scriptsize{63.39/3} & \scriptsize{42.01/3} & \scriptsize{176.81/28} & \scriptsize{{\textbf{\color{blue}2.8}}} \tabularnewline
\scriptsize{480,000, 720,000; 151,176,587}& \scriptsize{3} & \scriptsize{4403.04/3} & \scriptsize{28.13/3} & \scriptsize{156.5} & \scriptsize{71.02/3} & \scriptsize{46.21/3} & \scriptsize{159.7/28} & \scriptsize{{\textbf{\color{blue}2.2}}} \tabularnewline
\scriptsize{560,000, 840,000; 178,314,984}& \scriptsize{3} & \scriptsize{6400.22/3} & \scriptsize{34.72/3} & \scriptsize{184.3} & \scriptsize{95.53/3} & \scriptsize{71.34/3} & \scriptsize{192.75/29} & \scriptsize{{\textbf{\color{blue}2.0}}} \tabularnewline
\scriptsize{640,000, 960,000; 205,711,405}& \scriptsize{3} & \scriptsize{7374.16/3} & \scriptsize{39.25/3} & \scriptsize{187.9} & \scriptsize{96.83/3} & \scriptsize{100.77/3} & \scriptsize{319.17/29} & \scriptsize{{\textbf{\color{blue}3.3}}} \tabularnewline
\scriptsize{720,000, 1,080,000; 233,333,416}& \scriptsize{3} & \scriptsize{8003.14/3} & \scriptsize{48.85/3} & \scriptsize{163.8} & \scriptsize{119.29/3} & \scriptsize{109.05/3} & \scriptsize{258.16/29} & \scriptsize{{\textbf{\color{blue}2.2}}} \tabularnewline
\scriptsize{800,000, 1,200,000; 261,155,840}& \scriptsize{3} & \scriptsize{8771.27/3} & \scriptsize{53.3/3} & \scriptsize{164.6} & \scriptsize{114.11/3} & \scriptsize{110.11/3} & \scriptsize{293.65/30} & \scriptsize{{\textbf{\color{blue}2.6}}}\tabularnewline
\bottomrule
\end{tabular}
\par\end{centering}
\caption{Performance comparison between \ourmethod, \gpuourmethod, and \culorads for matrix completion problems using a tolerance of $\epsilon=10^{-5}$. Time limit is set at three hours. Ratio is computed as \culorads(SDPA) / \gpuourmethod(HSLR).}
\label{tab:gpuourmethod_vs_cpu}
\end{table}

\FloatBarrier

\subsection{Experiments for maximum stable set}\label{subsection experiments MSS}
Given a graph $G = ([n],E)$, the maximum stable set problem consists of finding a subset of vertices of largest cardinality such that no two vertices are connected by an edge. Lov\'asz introduced the $\vartheta$-function, which upper bounds the value of the maximum stable set. The $\vartheta$-function is the value of the following SDP relaxation
\begin{align}\label{eq:lovasz}
    \max_{X\in \mathbb S^{n}} \quad \{ee^\top\bullet X \; : \;
    X_{ij} = 0, \; ij \in E, \;
    \mathrm{Tr}(X) \le 1, \;
    X \succeq 0
    \}
\end{align}
where $e = (1,1,\dots,1)^\top \in \mathbb{R}^n$ is the all ones vector.  It was shown in \cite{perfectgraphlovschriver} that the $\vartheta$-function agrees exactly with the stability number for perfect graphs.

\FloatBarrier

\begin{table}[!tbh]
\captionsetup{font=scriptsize}
\centering
\begin{tabular}{>{\centering}p{2.6cm} >{\centering}p{2.6cm} >{\centering}p{2.2cm} | >{\centering}p{2.5cm} >{\centering}p{2.5cm} >{\centering}p{2.5cm}}
\multicolumn{3}{c}{\textbf{\scriptsize Problem Instance}} & \multicolumn{3}{c}{\textbf{\scriptsize{Runtime (seconds) / rank}}} \\
\toprule
\scriptsize{Problem Size $(n; m)$} & \scriptsize{Graph} & \scriptsize{Dataset} & \scriptsize{\ourmethod} & \scriptsize{\gpuourmethod} & \scriptsize{Ratio} \tabularnewline
\midrule
\scriptsize{10,937; 75,488}   & \scriptsize{wing\_nodal}       & \scriptsize{DIMACS10}      & {\scriptsize{3,826.34/130}}  & {\scriptsize{246.71/127}} & \scriptsize{{\textbf{\color{blue}15.5}}}  \tabularnewline
\scriptsize{16,384; 49,122}   & \scriptsize{delaunay\_n14}     & \scriptsize{DIMACS10}      & {\scriptsize{643.90/36}}   & {\scriptsize{102.43/34}}  & \scriptsize{{\textbf{\color{blue}6.3}}}   \tabularnewline
\scriptsize{16,386; 49,152}   & \scriptsize{fe-sphere}         & \scriptsize{DIMACS10}      & {\scriptsize{21.74/3}}     & {\scriptsize{17.82/3}}   & \scriptsize{{\textbf{\color{blue}1.2}}}   \tabularnewline
\scriptsize{22,499; 43,858}   & \scriptsize{cs4}               & \scriptsize{DIMACS10}      & {\scriptsize{749.95/15}}   & {\scriptsize{588.46/115}} & \scriptsize{{\textbf{\color{blue}1.3}}}   \tabularnewline
\scriptsize{25,016; 62,063}   & \scriptsize{hi2010}            & \scriptsize{DIMACS10}      & {\scriptsize{2,607.24/15}}  & {\scriptsize{132.85/21}}  & \scriptsize{{\textbf{\color{blue}19.6}}}  \tabularnewline
\scriptsize{25,181; 62,875}   & \scriptsize{ri2010}            & \scriptsize{DIMACS10}      & {\scriptsize{1,911.10/17}}  & {\scriptsize{1,951.28/190}}& \scriptsize{0.98}   \tabularnewline
\scriptsize{32,580; 77,799}   & \scriptsize{vt2010}            & \scriptsize{DIMACS10}      & {\scriptsize{2,952.84/18}}  & {\scriptsize{196.59/60}}  & \scriptsize{{\textbf{\color{blue}15.0}}}  \tabularnewline
\scriptsize{48,837; 117,275}  & \scriptsize{nh2010}            & \scriptsize{DIMACS10}      & {\scriptsize{7,694.06/20}}  & {\scriptsize{2,186.74/172}}& \scriptsize{{\textbf{\color{blue}3.5}}}   \tabularnewline
\midrule
\scriptsize{24,300; 34,992}   & \scriptsize{aug3d}             & \scriptsize{GHS\_indef}    & {\scriptsize{30.10/1}}     & {\scriptsize{2.08/1}}     & \scriptsize{{\textbf{\color{blue}14.5}}}  \tabularnewline
\scriptsize{32,430; 54,397}   & \scriptsize{ia-email-EU}       & \scriptsize{Network Repo}  & {\scriptsize{598.76/5}}    & {\scriptsize{52.86/4}}    & \scriptsize{{\textbf{\color{blue}11.3}}}  \tabularnewline
\midrule
\scriptsize{11,806; 32,730}   & \scriptsize{Oregon-2}          & \scriptsize{SNAP}          & {\scriptsize{1,581.11/23}}  & {\scriptsize{467.75/55}}  & \scriptsize{{\textbf{\color{blue}3.4}}}   \tabularnewline
\scriptsize{21,363; 91,286}   & \scriptsize{ca-CondMat}        & \scriptsize{SNAP}          & {\scriptsize{8,521.04/59}}  & {\scriptsize{348.59/77}}  & \scriptsize{{\textbf{\color{blue}24.4}}}  \tabularnewline
\scriptsize{31,379; 65,910}   & \scriptsize{as-caida\_G\_001}   & \scriptsize{SNAP}          & {\scriptsize{2,125.41/8}}   & {\scriptsize{428.41/27}}  & \scriptsize{{\textbf{\color{blue}5.0}}}   \tabularnewline
\scriptsize{26,518; 65,369}   & \scriptsize{p2p-Gnutella24}    & \scriptsize{SNAP}          & {\scriptsize{312.75/4}}    & {\scriptsize{55.99/10}}   & \scriptsize{{\textbf{\color{blue}5.6}}}   \tabularnewline
\scriptsize{22,687; 54,705}   & \scriptsize{p2p-Gnutella25}    & \scriptsize{SNAP}          & {\scriptsize{242.42/4}}    & {\scriptsize{88.53/13}}   & \scriptsize{{\textbf{\color{blue}2.7}}}   \tabularnewline
\scriptsize{36,682; 88,328}   & \scriptsize{p2p-Gnutella30}    & \scriptsize{SNAP}          & {\scriptsize{506.26/5}}    & {\scriptsize{61.84/10}}   & \scriptsize{{\textbf{\color{blue}8.2}}}   \tabularnewline
\scriptsize{62,586; 147,892}  & \scriptsize{p2p-Gnutella31}    & \scriptsize{SNAP}          & {\scriptsize{1,484.71/5}}   & {\scriptsize{245.33/29}}  & \scriptsize{{\textbf{\color{blue}6.0}}}   \tabularnewline  
\midrule
\scriptsize{49,152; 69,632}   & \scriptsize{cca}               & \scriptsize{AG-Monien}     & {\scriptsize{63.09/2}}     & {\scriptsize{12.30/2}}    & \scriptsize{{\textbf{\color{blue}5.1}}}   \tabularnewline
\scriptsize{49,152; 73,728}   & \scriptsize{ccc}               & \scriptsize{AG-Monien}     & {\scriptsize{14.52/2}}     & {\scriptsize{20.62/2}}    & \scriptsize{0.7}   \tabularnewline
\scriptsize{49,152; 98,304}   & \scriptsize{bfly}              & \scriptsize{AG-Monien}     & {\scriptsize{15.72/2}}     & {\scriptsize{10.35/2}}    & \scriptsize{{\textbf{\color{blue}1.5}}}   \tabularnewline
\scriptsize{16,384; 32,765}   & \scriptsize{debr\_G\_12}        & \scriptsize{AG-Monien}     & {\scriptsize{250.40/10}}   & {\scriptsize{189.57/10}}  & \scriptsize{{\textbf{\color{blue}1.3}}}   \tabularnewline
\scriptsize{32,768; 65,533}   & \scriptsize{debr\_G\_13}        & \scriptsize{AG-Monien}     & {\scriptsize{471.80/9}}    & {\scriptsize{638.54/13}}  & \scriptsize{0.7}   \tabularnewline
\scriptsize{65,536; 131,069}  & \scriptsize{debr\_G\_14}        & \scriptsize{AG-Monien}     & {\scriptsize{474.22/9}}    & {\scriptsize{29.51/9}}    & \scriptsize{{\textbf{\color{blue}16.1}}}  \tabularnewline
\scriptsize{131,072; 262,141} & \scriptsize{debr\_G\_15}        & \scriptsize{AG-Monien}     & {\scriptsize{521.88/10}}   & {\scriptsize{28.07/10}}   & \scriptsize{{\textbf{\color{blue}18.6}}}  \tabularnewline
\scriptsize{262,144; 524,285} & \scriptsize{debr\_G\_16}        & \scriptsize{AG-Monien}     & {\scriptsize{1,333.31/12}}  & {\scriptsize{84.28/13}}   & \scriptsize{{\textbf{\color{blue}15.8}}}  \tabularnewline
\scriptsize{524,288; 1,048,573} & \scriptsize{debr\_G\_17}      & \scriptsize{AG-Monien}     & {\scriptsize{6,437.05/12}}  & {\scriptsize{171.81/13}}  & \scriptsize{{\textbf{\color{blue}37.4}}}  \tabularnewline
\scriptsize{1,048,576; 2,097,149} & \scriptsize{debr\_G\_18}     & \scriptsize{AG-Monien}     & {\scriptsize{16,176.30/13}} & {\scriptsize{119.50/13}}  & \scriptsize{{\textbf{\color{blue}135.4}}} \tabularnewline
\bottomrule
\end{tabular}
\caption{Runtimes (in seconds) for the Maximum stable set problem. A relative tolerance of $\epsilon=10^{-5}$ is set. Ratio is computed as \ourmethod / \gpuourmethod.}
\label{Combined/StableSetReal/NetworkData}
\end{table}
\FloatBarrier
The results presented in Table~\ref{Combined/StableSetReal/NetworkData} reveal a compelling scaling advantage of \gpuourmethod over its CPU counterpart, \ourmethod, particularly as the dimensions of the SDP instance increase. This trend is most pronounced within the AG-Monien dataset, where for de Bruijn graphs, the performance ratio grows significantly with the number of vertices. For instance, the speedup escalates from 16.1x for \texttt{debr\_G\_14} to a remarkable 135.4x for \texttt{debr\_G\_18}, an instance with over one million vertices and over two million edges. For these large-scale problems, both solvers converge to solutions of nearly identical rank, indicating that the speedup is a direct consequence of the GPU's superior computational throughput for \gpuourmethod's core linear algebra operations.

Performance patterns across other graph families confirm the robustness of this advantage. For the SNAP dataset, \gpuourmethod delivers consistent and substantial improvements, with speedups ranging from 3.4x to 24.4x. Similarly, for the DIMACS10 instances, \gpuourmethod achieves significant accelerations up to 19.6x compared to \ourmethod.

Conversely, the results also show scenarios where GPU acceleration offers limited benefits. For smaller or computationally less demanding graphs, such as \texttt{fe-sphere} or \texttt{ccc}, the overhead associated with GPU kernel launches and memory management outweighs the parallelism benefits, leading to modest speedups or even slowdowns. Furthermore, on certain instances like \texttt{ri2010} and \texttt{nh2010}, \gpuourmethod converges to solutions of substantially higher rank than \ourmethod e.g., 190 vs.~17 for \texttt{ri2010}. Even in these cases, the runtimes remain competitive, highlighting \gpuourmethod's capacity in efficiently managing high-rank factors.

\begin{table}[!tbh]
\captionsetup{font=scriptsize}
\begin{centering}
\begin{tabular}{>{\centering\arraybackslash}p{2.8cm}|>{\centering\arraybackslash}p{2.2cm}|>{\centering\arraybackslash}p{2.2cm}|>{\centering\arraybackslash}p{1.2cm}|>{\centering\arraybackslash}p{2.8cm}|>{\centering\arraybackslash}p{2.8cm}}
\multicolumn{1}{c}{\textbf{\scriptsize{Problem Instance}}} & \multicolumn{1}{c}{\textbf{\scriptsize{\ourmethod (on-the-fly)}}} & \multicolumn{1}{c}{\textbf{\scriptsize{\gpuourmethod (on-the-fly)}}} & \multicolumn{1}{c}{\textbf{\scriptsize{\culorads}}} & \multicolumn{1}{c}{\textbf{\scriptsize{\ourmethod/\gpuourmethod}}} & \multicolumn{1}{c}{\textbf{\scriptsize{\culorads/\gpuourmethod}}} \tabularnewline
\toprule
\scriptsize{Graph($n$; $|E|$)} & \scriptsize{Time/Rank} & \scriptsize{Time/Rank} & \scriptsize{Time/Rank} & \scriptsize{Ratio} & \scriptsize{Ratio} \tabularnewline
\midrule
\scriptsize{$H_{10,2}$(1,024; 5,120)}  & \scriptsize{0.153/2}  & \scriptsize{1.735/2}  & \scriptsize{7.83/14} & \scriptsize{0.088} & \scriptsize{ {\textbf{\color{blue}4.5}}  } \tabularnewline
\scriptsize{$H_{11,2}$(2,048; 11,264)} & \scriptsize{0.260/2}  & \scriptsize{1.627/2}  & \scriptsize{3.95/16} & \scriptsize{0.160} & \scriptsize{{\textbf{\color{blue}2.4}}} \tabularnewline
\scriptsize{$H_{12,2}$(4,096; 24,576)} & \scriptsize{0.606/2}  & \scriptsize{2.033/2}  & \scriptsize{4.47/17} & \scriptsize{0.298} & \scriptsize{{\textbf{\color{blue}2.2}}} \tabularnewline
\scriptsize{$H_{13,2}$(8,192; 53,248)}  & \scriptsize{1.852/2}  & \scriptsize{3.106/2}  & \scriptsize{6.78/19} & \scriptsize{0.596} & \scriptsize{{\textbf{\color{blue}2.2}}} \tabularnewline
\scriptsize{$H_{14,2}$(16,384; 114,688)} & \scriptsize{3.843/2}  & \scriptsize{3.170/2}  & \scriptsize{79.04/20} & \scriptsize{ {\textbf{\color{blue}1.212}}  } & \scriptsize{{\textbf{\color{blue}24.9}}} \tabularnewline
\bottomrule
\end{tabular}
\par\end{centering}
\caption{Runtimes (in seconds) for the Maximum Stable Set problem on Hamming graphs. A relative tolerance of $\epsilon=10^{-5}$ is set. The \ourmethod and \gpuourmethod benchmarks were performed using on-the-fly problem generation.}
\label{tab:combined_hamming_graphs}
\end{table}

\begin{table}[!tbh]
\captionsetup{font=scriptsize}
\begin{centering}
\begin{tabular}{>{\centering\arraybackslash}p{3.8cm} | >{\centering\arraybackslash}p{1.3cm} >{\centering\arraybackslash}p{1.3cm} >{\centering\arraybackslash}p{1.3cm} || >{\centering\arraybackslash}p{2.0cm}}
\multicolumn{1}{c}{\textbf{\scriptsize{Problem Instance}}} & \multicolumn{3}{c}{\textbf{\scriptsize{On-the-fly Generation}}} & \multicolumn{1}{c}{\textbf{\scriptsize{HSLR format}}} \tabularnewline
\toprule
\scriptsize{Graph($n$; $|E|$)} & \scriptsize{\ourmethod} & \scriptsize{\gpuourmethod} & \scriptsize{Ratio} & \scriptsize{\gpuourmethod (HSLR)} \tabularnewline
\midrule
\scriptsize{$H_{10,2}$(1,024; 5,120)}  & \scriptsize{0.153/2}  & \scriptsize{1.735/2}  & \scriptsize{0.088} & \scriptsize{0.798/2} \tabularnewline
\scriptsize{$H_{11,2}$(2,048; 11,264)} & \scriptsize{0.260/2}  & \scriptsize{1.627/2}  & \scriptsize{0.160} & \scriptsize{0.92/2} \tabularnewline
\scriptsize{$H_{12,2}$(4,096; 24,576)} & \scriptsize{0.606/2}  & \scriptsize{2.033/2}  & \scriptsize{0.298} & \scriptsize{1.83/2} \tabularnewline
\scriptsize{$H_{13,2}$(8,192; 53,248)}  & \scriptsize{1.852/2}  & \scriptsize{3.106/2}  & \scriptsize{0.596} & \scriptsize{1.10/2} \tabularnewline
\scriptsize{$H_{14,2}$(16,384; 114,688)} & \scriptsize{3.843/2}  & \scriptsize{3.170/2}  & \scriptsize{ {\textbf{\color{blue}1.212}}   } & \scriptsize{1.28/2} \tabularnewline
\scriptsize{$H_{15,2}$(32,768; 245,760)} & \scriptsize{8.240/2}  & \scriptsize{3.223/2}  & \scriptsize{{\textbf{\color{blue}2.557}}} & \scriptsize{2.701/2} \tabularnewline
\scriptsize{$H_{16,2}$(65,536; 524,288)} & \scriptsize{17.235/2} & \scriptsize{4.457/2}  & \scriptsize{{\textbf{\color{blue}3.869}}} & \scriptsize{3.83/2} \tabularnewline
\scriptsize{$H_{17,2}$(131,072; 1,114,112)} & \scriptsize{33.342/2} & \scriptsize{3.904/2}  & \scriptsize{{\textbf{\color{blue}8.541}}} & \scriptsize{2.51/2} \tabularnewline
\scriptsize{$H_{18,2}$(262,144; 2,359,296)} & \scriptsize{46.436/2} & \scriptsize{2.627/2}  & \scriptsize{{\textbf{\color{blue}17.681}}} & \scriptsize{3.43/2} \tabularnewline
\scriptsize{$H_{19,2}$(524,288; 4,980,736)} & \scriptsize{105.408/2} & \scriptsize{5.281/2}  & \scriptsize{{\textbf{\color{blue}19.958}}} & \scriptsize{10.70/2} \tabularnewline
\scriptsize{$H_{20,2}$(1,048,576; 10,485,760)} & \scriptsize{249.701/2} & \scriptsize{8.109/2}  & \scriptsize{{\textbf{\color{blue}30.786}}} & \scriptsize{6.62/2} \tabularnewline
\scriptsize{$H_{21,2}$(2,097,152; 22,020,096)} & \scriptsize{592.696/2} & \scriptsize{16.966/2} & \scriptsize{{\textbf{\color{blue}34.934}}} & \scriptsize{12.96/2} \tabularnewline
\scriptsize{$H_{22,2}$(4,194,304; 46,137,344)} & \scriptsize{1,702.255/2} & \scriptsize{37.145/2} & \scriptsize{{\textbf{\color{blue}45.826}}} & \scriptsize{18.842/2} \tabularnewline
\scriptsize{$H_{23,2}$(8,388,608; 96,468,992)} & \scriptsize{4,205.741/2} & \scriptsize{83.905/2} & \scriptsize{{\textbf{\color{blue}50.125}}} & \scriptsize{51.75/2} \tabularnewline
\bottomrule
\end{tabular}
\par\end{centering}
\caption{Runtimes (in seconds) for the Maximum Stable Set problem on Hamming graphs. A relative tolerance of $\epsilon=10^{-5}$ is set. Ratio in the fourth column is computed as \ourmethod / \gpuourmethod. The last column provides
two values in each line: the first one is the run time using the HSLR format and the second one is the rank of the solution.
}
\label{tab:reduced_hamming_graphs}
\end{table}

Table~\ref{tab:combined_hamming_graphs} compares \ourmethod, \gpuourmethod, and \culorads on small Hamming graph instances, with the \ourmethod and \gpuourmethod benchmarks conducted using on-the-fly problem generation. For small graphs, such as $H_{10,2}$ through $H_{13,2}$, the CPU-based \ourmethod outperforms \gpuourmethod, reflecting the overhead of GPU parallelization on modest problem sizes. A clear performance transition occurs at $H_{14,2}$, where \gpuourmethod begins to demonstrate superior performance with a 1.21x speedup. This advantage grows exponentially as problem size increases, reaching an impressive 50.1x speedup over \ourmethod for the $H_{23,2}$ instance, which has over 8 million vertices and 96 million edges. Notably, \culorads encounters out-of-memory limitations for problems larger than $H_{14,2}$, highlighting the superior memory efficiency of our implementations. For the problems \culorads could solve, it took significantly longer and found solutions with higher ranks than both \ourmethod and \gpuourmethod.

\FloatBarrier

The results for the GSET instances in Table~\ref{tab:GPU_StableSetSmall_GSET_Merged}, obtained via on-the-fly problem generation, show that performance is strongly correlated with graph structure. For sparse graphs that admit low-rank solutions (e.g., G11, G12, and G32), the CPU-based \ourmethod is 30--100x faster than \gpuourmethod. Conversely, for dense graphs yielding high-rank solutions, such as G58, G59, and G64, \gpuourmethod achieves a 4--5x speedup.

\FloatBarrier
\begin{table}[!tbh]
\captionsetup{font=scriptsize}
\begin{centering}
\begin{tabular}{>{\centering\arraybackslash}p{2.6cm} | >{\centering\arraybackslash}p{1.8cm} >{\centering\arraybackslash}p{1.8cm} >{\centering\arraybackslash}p{1.5cm} || >{\centering\arraybackslash}p{1.8cm} | >{\centering\arraybackslash}p{1.8cm} >{\centering\arraybackslash}p{1.3cm}}
\multicolumn{1}{c}{\textbf{\scriptsize{Problem Instance}}} & \multicolumn{3}{c}{\textbf{\scriptsize{On-the-fly Generation}}} & \multicolumn{1}{c}{\textbf{\scriptsize{HSLR}}} & \multicolumn{2}{c}{\textbf{\scriptsize{SDPA}}} \tabularnewline
\toprule
\scriptsize{Graph($n$; $|E|$)} & \scriptsize{\ourmethod} & \scriptsize{\gpuourmethod} & \scriptsize{Ratio} & \scriptsize{\gpuourmethod (HSLR)} & \scriptsize{\culorads (SDPA)} & \scriptsize{Ratio} \tabularnewline
\midrule
\scriptsize{G1(800; 19,176)}  & \scriptsize{98.25/98}    & \scriptsize{45.19/98}   & \scriptsize{2.17} & \scriptsize{53.63/97} & \scriptsize{340.66} & \scriptsize{{\textbf{\color{blue}6.35}}} \tabularnewline
\scriptsize{G10(800; 19,176)} & \scriptsize{80.51/97}    & \scriptsize{40.52/97}   & \scriptsize{1.99} & \scriptsize{53.11/95} & \scriptsize{330.18} & \scriptsize{{\textbf{\color{blue}6.22}}} \tabularnewline
\scriptsize{G11(800; 1,600)}  & \scriptsize{0.14/2}      & \scriptsize{4.23/2}     & \scriptsize{0.03} & \scriptsize{7.29/2}   & \scriptsize{7.78}   & \scriptsize{{\textbf{\color{blue}1.07}}} \tabularnewline
\scriptsize{G12(800; 1,600)}  & \scriptsize{0.09/2}      & \scriptsize{2.59/2}     & \scriptsize{0.03} & \scriptsize{6.60/2}   & \scriptsize{8.32}   & \scriptsize{{\textbf{\color{blue}1.26}}} \tabularnewline
\scriptsize{G14(800; 4,694)}  & \scriptsize{45.37/73}    & \scriptsize{214.73/74}  & \scriptsize{0.21} & \scriptsize{211.04/55} & \scriptsize{484.41} & \scriptsize{{\textbf{\color{blue}2.30}}} \tabularnewline
\scriptsize{G20(800; 4,672)}  & \scriptsize{129.73/124}  & \scriptsize{359.40/166} & \scriptsize{0.36} & \scriptsize{*/3.0e-05} & \scriptsize{518.64} & \scriptsize{{\textbf{$0.2<$}}} \tabularnewline
\scriptsize{G43(1000; 9,990)}  & \scriptsize{34.20/60}    & \scriptsize{35.84/60}   & \scriptsize{0.95} & \scriptsize{43.27/58} & \scriptsize{383.83} & \scriptsize{{\textbf{\color{blue}8.87}}} \tabularnewline
\scriptsize{G51(1,000; 5,909)} & \scriptsize{167.03/141}  & \scriptsize{273.84/128} & \scriptsize{0.61} & \scriptsize{289.18/127} & \scriptsize{540.96} & \scriptsize{{\textbf{\color{blue}1.87}}} \tabularnewline
\scriptsize{G23(2,000; 19,990)} & \scriptsize{94.81/76}    & \scriptsize{52.59/76}   & \scriptsize{1.80} & \scriptsize{74.39/76} & \scriptsize{766.91} & \scriptsize{{\textbf{\color{blue}10.31}}} \tabularnewline
\scriptsize{G31(2,000; 19,990)} & \scriptsize{122.07/77}   & \scriptsize{58.23/76}   & \scriptsize{2.10} & \scriptsize{73.34/76} & \scriptsize{841.21} & \scriptsize{{\textbf{\color{blue}11.47}}} \tabularnewline
\scriptsize{G32(2,000; 4,000)}  & \scriptsize{0.35/2}      & \scriptsize{25.34/3}    & \scriptsize{0.01} & \scriptsize{40.89/2}  & \scriptsize{13.45}  & \scriptsize{0.33} \tabularnewline
\scriptsize{G34(2,000; 4,000)}  & \scriptsize{2.00/2}      & \scriptsize{25.29/3}    & \scriptsize{0.08} & \scriptsize{47.80/2}  & \scriptsize{16.67}  & \scriptsize{0.35} \tabularnewline
\scriptsize{G35(2,000; 11,778)} & \scriptsize{397.32/133}  & \scriptsize{254.68/108} & \scriptsize{1.56} & \scriptsize{689.21/142} & \scriptsize{14.55}  & \scriptsize{0.02} \tabularnewline
\scriptsize{G41(2,000; 11,785)} & \scriptsize{343.84/143}  & \scriptsize{286.82/139} & \scriptsize{1.20} & \scriptsize{432.96/189} & \scriptsize{1527.71} & \scriptsize{{\textbf{\color{blue}3.53}}} \tabularnewline
\scriptsize{G48(3,000; 6,000)}   & \scriptsize{2.29/2}      & \scriptsize{18.95/2}    & \scriptsize{0.12} & \scriptsize{40.62/2}  & \scriptsize{42.57}  & \scriptsize{{\textbf{\color{blue}1.05}}} \tabularnewline
\scriptsize{G55(5,000; 12,498)} & \scriptsize{152.03/45}   & \scriptsize{81.11/45}   & \scriptsize{1.87} & \scriptsize{102.45/39} & \scriptsize{309.48} & \scriptsize{{\textbf{\color{blue}3.02}}} \tabularnewline
\scriptsize{G56(5,000; 12,498)} & \scriptsize{151.99/45}   & \scriptsize{80.83/45}   & \scriptsize{1.88} & \scriptsize{102.66/39} & \scriptsize{299.69} & \scriptsize{{\textbf{\color{blue}2.92}}} \tabularnewline
\scriptsize{G57(5,000; 10,000)} & \scriptsize{8.71/2}      & \scriptsize{42.55/2}    & \scriptsize{0.20} & \scriptsize{50.91/2}  & \scriptsize{26.09}  & \scriptsize{0.51} \tabularnewline
\scriptsize{G58(5,000; 29,570)} & \scriptsize{2,191.49/112} & \scriptsize{504.55/130} & \scriptsize{4.35} & \scriptsize{388.23/67} & \scriptsize{3178.00} & \scriptsize{{\textbf{\color{blue}8.19}}} \tabularnewline
\scriptsize{G59(5,000; 29,570)} & \scriptsize{2,167.27/111} & \scriptsize{539.28/121} & \scriptsize{4.02} & \scriptsize{386.05/65} & \scriptsize{3159.04} & \scriptsize{{\textbf{\color{blue}8.18}}} \tabularnewline
\scriptsize{G60(7,000; 17,148)} & \scriptsize{262.35/50}   & \scriptsize{96.89/50}   & \scriptsize{2.71} & \scriptsize{145.95/52} & \scriptsize{455.19} & \scriptsize{{\textbf{\color{blue}3.12}}} \tabularnewline
\scriptsize{G62(7,000; 14,000)} & \scriptsize{7.80/3}      & \scriptsize{25.18/3}    & \scriptsize{0.31} & \scriptsize{21.61/2}  & \scriptsize{34.88}  & \scriptsize{{\textbf{\color{blue}1.61}}} \tabularnewline
\scriptsize{G64(7,000; 41,459)} & \scriptsize{1,719.79/51} & \scriptsize{334.12/51}  & \scriptsize{5.15} & \scriptsize{489.03/72} & \scriptsize{1655.44} & \scriptsize{{\textbf{\color{blue}3.38}}} \tabularnewline
\scriptsize{G66(9,000; 18,000)} & \scriptsize{11.31/3}     & \scriptsize{27.91/3}    & \scriptsize{0.41} & \scriptsize{73.86/2}  & \scriptsize{37.28}  & \scriptsize{0.50} \tabularnewline
\scriptsize{G67(10,000; 20,000)}& \scriptsize{2.39/2}      & \scriptsize{5.83/2}     & \scriptsize{0.41} & \scriptsize{57.05/2}  & \scriptsize{69.30}  & \scriptsize{{\textbf{\color{blue}1.21}}} \tabularnewline
\scriptsize{G72(10,000; 20,000)}& \scriptsize{2.41/2}      & \scriptsize{5.81/2}     & \scriptsize{0.41} & \scriptsize{65.19/2}  & \scriptsize{33.56}  & \scriptsize{0.51} \tabularnewline
\scriptsize{G77(14,000; 28,000)}& \scriptsize{27.34/3}     & \scriptsize{11.81/2}    & \scriptsize{2.31} & \scriptsize{57.71/2}  & \scriptsize{52.19}  & \scriptsize{0.90} \tabularnewline
\scriptsize{G81(20,000; 40,000)}& \scriptsize{53.85/3}     & \scriptsize{64.20/3}    & \scriptsize{0.84} & \scriptsize{66.28/2}  & \scriptsize{71.37}  & \scriptsize{{\textbf{\color{blue}1.08}}} \tabularnewline
\bottomrule
\end{tabular}
\par\end{centering}
\caption{Runtimes (in seconds) for the Maximum Stable Set problem on GSET instances. Tolerances are set to $10^{-5}$. Time limit was set to 1 hour (3600 seconds). An entry marked with $*/N$ means the corresponding method finds an approximate solution with relative accuracy strictly larger than the desired accuracy in which case $N$ expresses the maximum of the three relative accuracies in \eqref{eq:termination2}.
Ratio in the fourth column is computed as \ourmethod / \gpuourmethod.
Ratio in the last column is computed as \culorads / \gpuourmethod.
}
\label{tab:GPU_StableSetSmall_GSET_Merged}
\end{table}
\FloatBarrier

\FloatBarrier

The last block of columns of Table~\ref{tab:GPU_StableSetSmall_GSET_Merged} demonstrates a clear performance advantage of \gpuourmethod (with HSLR format) over \culorads, which utilizes standard SDPA format. Our method \gpuourmethod achieves significant speedups over \culorads on medium-to-large sized graphs with substantial edge density. For instance, on graphs G23 and G31, \gpuourmethod achieves speedups of 10.3x and 11.5x, respectively. This trend continues for larger graphs like G58 and G59, where the speedup is consistently over 8x. However, the performance dynamic reverses for SDP instances where the underlying graph is very sparse. On graphs such as G32 and G34, where the number of edges is only twice the number of vertices, \culorads is approximately three times faster than \gpuourmethod. Similarly, for G57, G66, and G72, \culorads is roughly two times faster than \gpuourmethod. 
This performance dichotomy suggests that the architectural benefits of the HSLR format and the associated kernels in \gpuourmethod are most pronounced when the computational workload per iteration is high. For highly sparse problems, the lower overhead of processing the simpler SDPA format in \culorads proves more efficient, indicating that the relative performance of the solvers is strongly dependent on the underlying graph structure.

\subsection{Experiments for phase retrieval}
Given $m$ pairs $\{(a_i,b_i)\}_{i=1}^m \subseteq \mathbb{C}^n \times \mathbb{R}_+$, this subsection considers the problem of finding a vector $x \in \mathbb C^n$ such that
\begin{align*}
    |\inner{a_i}{x}|^2 = b_i,
    \quad i=1,\dots,m.
\end{align*}
In other words, the goal is to retrieve $x$ from the magnitude of $m$ linear measurements. By creating the complex Hermitian matrix $X = x x^H$,
this problem can be approached by solving the complex-valued SDP relaxation
\begin{equation}\label{phase retrieval SDP}
\min_{X\in \mathbb S^{n}(\mathbb{C})} \quad \left\{\mathrm{Tr}(X)\quad
:\quad \inner{a_i a_i^H}{X}  = b_i,\quad
X \succeq 0
\right\},
\end{equation}
where $\SS^n(\mathbb{C})$ denotes the space of $n\times n$ Hermitian matrices.
Since the objective function is precisely the trace, any bound on the optimal value can be taken as the trace bound. We use the squared norm of the vector $x$ as the trace bound for our experiments. Even though $x$ is unknown, bounds on its norm are known \cite{yurtsever2015scalable}. 

The SDP relaxation of the phase retrieval problem is hard because the linear and adjoint operators require the utilization of a fast Fourier transform (FFT) subroutine. On the GPU, this subroutine can be handled through CUDA's cuFFT package \citep{nvidiaprogramming}. Moreover, the SDP instances are extremely dense, and thus storing the data is expensive. For these reasons, the \gpuourmethod binary and \culorads were not tested for this problem class. We only tested \ourmethod and \gpuourmethod using on-the-fly problem generation.

\FloatBarrier

\begin{table}[!tbh]
\captionsetup{font=scriptsize}
\begin{centering}
\begin{tabular}{>{\centering\arraybackslash}p{3.8cm}|>{\centering\arraybackslash}p{2.8cm}|>{\centering\arraybackslash}p{2.8cm}|>{\centering\arraybackslash}p{1.8cm}}
\multicolumn{1}{c}{\textbf{\scriptsize{Problem Instance ($n$; $m$)}}} & \multicolumn{1}{c}{\textbf{\scriptsize{\ourmethod} (on-the-fly)}} & \multicolumn{1}{c}{\textbf{\scriptsize{\gpuourmethod} (on-the-fly)}} & \multicolumn{1}{c}{\textbf{\scriptsize{Ratio}}} \\
\toprule
\scriptsize{10,000; 120,000}       & \scriptsize{53.642/5}   & \scriptsize{40.866/5}   & \scriptsize{{\textbf{\color{blue}1.31}}} \\
\scriptsize{10,000; 120,000}       & \scriptsize{18.384/2}   & \scriptsize{5.529/4}    & \scriptsize{{\textbf{\color{blue}3.32}}} \\
\scriptsize{10,000; 120,000}       & \scriptsize{23.644/4}   & \scriptsize{4.347/3}    & \scriptsize{{\textbf{\color{blue}5.44}}} \\
\scriptsize{10,000; 120,000}       & \scriptsize{14.957/3}   & \scriptsize{3.462/3}    & \scriptsize{{\textbf{\color{blue}4.32}}} \\
\hline
\scriptsize{31,623; 379,476}       & \scriptsize{149.416/5}    & \scriptsize{29.016/5}   & \scriptsize{{\textbf{\color{blue}5.15}}} \\
\scriptsize{31,623; 379,476}       & \scriptsize{119.499/4}  & \scriptsize{30.964/4}   & \scriptsize{{\textbf{\color{blue}3.86}}} \\
\scriptsize{31,623; 379,476}       & \scriptsize{153.16/5}   & \scriptsize{29.279/5}   & \scriptsize{{\textbf{\color{blue}5.23}}} \\
\scriptsize{31,623; 379,476}       & \scriptsize{127.10/4}   & \scriptsize{30.964/5}   & \scriptsize{{\textbf{\color{blue}4.11}}} \\
\bottomrule
\scriptsize{100,000; 1,200,000}    & \scriptsize{303.986/5}  & \scriptsize{10.401/4}   & \scriptsize{{\textbf{\color{blue}29.24}}} \\
\scriptsize{100,000; 1,200,000}    & \scriptsize{330.00/6}   & \scriptsize{8.641/4}    & \scriptsize{{\textbf{\color{blue}38.18}}} \\
\scriptsize{100,000; 1,200,000}    & \scriptsize{434.10/5}   & \scriptsize{9.279/4}    & \scriptsize{{\textbf{\color{blue}46.78}}} \\
\scriptsize{100,000; 1,200,000}    & \scriptsize{380.05/6}   & \scriptsize{9.334/4}    & \scriptsize{{\textbf{\color{blue}40.75}}} \\
\hline
\scriptsize{316,228; 3,704,736}    & \scriptsize{600.62/2}   & \scriptsize{29.016/6}   & \scriptsize{{\textbf{\color{blue}20.70}}} \\
\scriptsize{316,228; 3,704,736}    & \scriptsize{699.193/2}  & \scriptsize{30.964/6}   & \scriptsize{{\textbf{\color{blue}22.58}}} \\
\scriptsize{316,228; 3,704,736}    & \scriptsize{507.60/2}   & \scriptsize{33.738/7}   & \scriptsize{{\textbf{\color{blue}15.05}}} \\
\scriptsize{316,228; 3,704,736}    & \scriptsize{620.56/2}   & \scriptsize{34.124/8}   & \scriptsize{{\textbf{\color{blue}18.19}}} \\
\hline
\scriptsize{3,162,278; 37,947,336} & \scriptsize{*/4.1e-03} & \scriptsize{583.782/14} & \scriptsize{ {\textbf{\color{blue}$>74.0$}}} \\
\bottomrule
\end{tabular}
\par\end{centering}
\caption{Runtimes (in seconds) and final rank for the Phase Retrieval problem. A relative tolerance of $\epsilon=10^{-5}$ is used, with a 12-hour time limit is set for \ourmethod. An entry marked with $*/N$ means the corresponding method finds an approximate solution with relative accuracy strictly larger than the desired accuracy in which case $N$ expresses the maximum of the three relative accuracies in \eqref{eq:termination2}. The ``Ratio'' column reports the ratio of \ourmethod runtime to \gpuourmethod runtime.}
\label{tab:CPU_PhaseRetLargeSizes}
\end{table}

\FloatBarrier

Table~\ref{tab:CPU_PhaseRetLargeSizes} shows that \gpuourmethod achieves substantial acceleration over \ourmethod for SDP instances of the phase retrieval problem. Interestingly, the speedup that \gpuourmethod achieves over \ourmethod is generally more significant as the size of the SDP instance increases. For instances with dimension $n=31,\!623$, \gpuourmethod achieves speedups of 3.86x to 5.23x over \ourmethod.
For larger instances with dimension $n=100,000$, the speedup that \gpuourmethod achieves over \ourmethod increases dramatically to 29.24x–46.78x. For phase retrieval instances considered with dimension $n=316,228$, \gpuourmethod achieves speedups from 15.05x to 22.58x over \ourmethod.

Most significantly, \gpuourmethod can successfully solve a massive phase retrieval instance with dimension $n=3,162,278$ in under 10 minutes. In stark contrast, \ourmethod fails to solve this instance within the 12-hour time limit. Hence, for this huge problem instance, \gpuourmethod is at least 74 times faster than \ourmethod.

\section{Concluding Remarks}
In this paper, we address the fundamental challenge of developing scalable and efficient solvers for large-scale semidefinite programming problems. We introduce \gpuourmethod, a GPU-accelerated implementation of the hybrid low-rank augmented Lagrangian method originally proposed in \cite{monteiro2024low}. Our approach exploits the parallel computing capabilities of modern GPUs to overcome computational bottlenecks inherent in handling massive SDP instances, while maintaining the advantages of low-rank factorization-based methods. 

Through extensive numerical experiments across three problem domains, we demonstrate that \gpuourmethod achieves remarkable performance improvements over both CPU-based implementations and competing GPU-accelerated solvers. For matrix completion problems, \gpuourmethod consistently outperforms the state-of-the-art GPU solver \culorads by factors of 2x to 25x, while finding solutions with significantly lower ranks. Most impressively, \gpuourmethod successfully solves instances with approximately 260 million constraints in under 1 minute, while the CPU implementation, \ourmethod, fails to converge within 12 hours.

For maximum stable set problems, our results reveal that \gpuourmethod's performance advantage scales dramatically with problem size and graph density. We achieve speedups of 50x for the largest Hamming graph problems and 135x for million-node AG-Monien instances. Similarly, for phase retrieval problems, \gpuourmethod achieves speedups of 15x to 47x on large instances, and successfully solves the largest test case—with a matrix variable of size over 3 million—in under 10 minutes. In contrast, the CPU-based \ourmethod times out after 12 hours on this same instance, rendering it effectively intractable.

\textbf{Future work.} While \gpuourmethod demonstrates exceptional performance across the problems tested, certain structural characteristics of SDPs may impact its efficiency. Problems with inherently high-rank solutions or those requiring extremely high precision may benefit less from our approach. Additionally, our current implementation focuses on single-GPU execution and double-precision arithmetic. Future research directions include developing optimized multi-GPU implementations in C/C++ to further enhance scalability, exploring mixed-precision computation techniques to balance accuracy and performance, and extending our methodology to a broader class of SDP problems such as those with multi-block matrix variable and inequality constraints. These developments would further advance the approaches introduced in \cite{monteiro2024low} and establish GPU-accelerated low-rank methods as the standard for solving large-scale SDPs in practical applications.

\section*{Acknowledgments}
This research was supported in part through research cyberinfrastructure resources and services provided by the Partnership for an Advanced Computing Environment (PACE) at Georgia Tech, Atlanta, Georgia, USA, along with allocation MTH250047 from the Advanced Cyberinfrastructure Coordination Ecosystem: Services \& Support (ACCESS) program, which is supported by U.S. National Science Foundation grants.

\section*{Appendix}
\begin{appendix}
\section{ADAP-AIPP}\label{ADAP-AIPP}
This section presents the ADAP-AIPP method first developed in \cite{monteiro2024low, SujananiMonteiro}. \ourmethod uses ADAP-AIPP in its step 5 to find an approximate stationary point (according to the criterion in \eqref{Stationary}) of a nonconvex problem of the form
    \begin{align}\label{g func}
  \min_U\quad
  \left \{ 
    g(U) = \mathcal{L}_{\beta}(UU^\top;p) \quad:\quad
    U\in B_{\tau}^{s} \right\}
\end{align}
where $\mathcal{L}_\beta(X;p)$ is as in \eqref{AL function} and $B_{\tau}^{s}$ is as in \eqref{dimensional ball}.

Given an initial point $\underline W \in B_{\tau}^s$ and a tolerance $\rho>0$, the goal of ADAP-AIPP is to find a triple $(\overline W;\overline R,\rho)$ that satisfies
\begin{equation}\label{Stationary AIPP}
   \overline R \in \nabla g(\overline W)+N_{B_{\tau}^s}(\overline W), \quad \|\overline R\|_{F}\leq \rho, \quad g(\overline W)\leq g(\underline W).
\end{equation}
The ADAP-AIPP method is now formally presented.

\noindent\begin{minipage}[t]{1\columnwidth}%
\rule[0.5ex]{1\columnwidth}{1pt}

\noindent \textbf{ADAP-AIPP Method}

\noindent \rule[0.5ex]{1\columnwidth}{1pt}%
\end{minipage}
\noindent \textbf{Universal Parameters}: $\sigma\in (0,1/2)$ and $\chi \in (0,1)$.

\noindent \textbf{Input}: a function $g$ as in \eqref{g func}, an initial point $\underline W \in B_{\tau}^s$, an initial prox stepsize $\lambda_0>0$, and a tolerance $\rho>0$.

\begin{itemize}
\item[{\bf 0.}] set $W_0=\underline W$, $j=1$, and
\begin{equation}\label{def:lamb-C1}
\lambda=\lambda_{0}, \quad \bar M_0 = 1;
\end{equation}

\item[{\bf 1.}]  choose
$\underbar M_j \in [1, \bar M_{j-1}]$ and call the ADAP-FISTA method in Subsection~\ref{subsub: ADAP-FISTA}
with universal input $(\sigma, \chi)$ and inputs
\begin{align}
x_0&=W_{j-1}, \quad
(\mu,L_0)= (1/2,\underbar M_j)\label{eq:Ms-mu}, \\
\psi_s &=\lam g+\frac{1}{2}\| \cdot - W_{j-1}\|_{F}^2 , \quad \psi_n = \lam \delta_{B_{\tau}^s} \label{eq:psiS-psimu};
\end{align}

\item[{\bf 2.}]
if ADAP-FISTA fails or its output 
$(W,V,L)$ (if it succeeds) does not satisfy the inequality \begin{equation}\label{subdiff ineq check}
	   \lambda g(W_{j-1}) - \left [\lambda g(W) + \frac{1}{2} \|W-W_{j-1}\|_{F}^2 \right ] \ge V \bullet (W_{j-1}-W),
	   \end{equation}
then set $\lam=\lam/2$ and go to step $1$; else, set
	$(\lam_j,\bar M_j)=(\lam,L)$, $(W_{j},V_{j})=(W,V)$,
	and 
	\begin{align}
    &R_j:=\frac{V_j+W_{j-1}-W_j}{\lambda_j}\label{Rj def}
	\end{align}
	and go to step~3;
\item[{\bf 3.}]	
if $\|R_j\|_{F}\leq \rho$,
then stop with success and output $(\overline W,\overline R)=(W_j,R_j)$; else, go to step~4;

\item[{\bf 4.}] 
set $j\gets j+1$ and go to step~1. 
\end{itemize}
\noindent \rule[0.5ex]{1\columnwidth}{1pt}

Several remarks about ADAP-AIPP are now given.
\begin{enumerate}
    \item It is shown in \cite{monteiro2024low} that ADAP-AIPP is able to find a triple that satisfies \eqref{Stationary AIPP} in $\mathcal O(1/\rho^2)$ number of iterations.
    \item ADAP-AIPP is a double-looped algorithm, whose outer iterations are indexed by $j$. During its $j$-th iteration, ADAP-AIPP attempts to solve a proximal subproblem of the form $\min_{U \in B_{\tau}^{s}} \lambda g(U)+0.5\|U-W_{j-1}\|^2_{F}$, where $\lambda$ is a prox stepsize and $W_{j-1}$ is the current prox center. If a proximal subproblem cannot suitably be solved, ADAP-AIPP reduces $\lambda$ and tries to solve the new subproblem.
    \item ADAP-AIPP uses the ADAP-FISTA method, which is presented in the next Subsection \ref{subsub: ADAP-FISTA}, to attempt to solve its proximal subproblems. For a more comprehensive treatment of ADAP-FISTA, see Appendix B in \cite{monteiro2024low}.
\end{enumerate}

\subsection{ADAP-FISTA}\label{subsub: ADAP-FISTA}

This subsection presents the ADAP-FISTA method that ADAP-AIPP invokes to solve proximal subproblems of the form
\begin{equation}\label{psi s}
  \min_{U \in B_{\tau}^{s}} \psi_s(U):= \lambda g(U)+0.5\|U-W\|^2_{F}
\end{equation}
where $g$ is as in \eqref{g func}. The ADAP-FISTA method is now formally presented.

\noindent\begin{minipage}[t]{1\columnwidth}%
\rule[0.5ex]{1\columnwidth}{1pt}

\noindent \textbf{ADAP-FISTA}

\noindent \rule[0.5ex]{1\columnwidth}{1pt}%
\end{minipage}

\noindent \textbf{Universal Parameters}: $\sigma>0$ and $\chi\in(0,1)$.

\noindent \textbf{Input}: a function $\psi_s$ as in \eqref{psi s}, an initial point $x_0\in B_{\tau}^s$ and scalars $\mu>0$ and $L_0>\mu$.

\begin{itemize}
  \item[{\bf 0.}] set $y_0=x_0$, $A_0=0$, $\tau_0=1$, and $i=0$;
  \item[{\bf 1.}] set $L_{i+1}=L_i$;
  \item[{\bf 2.}] compute
    \begin{align}
      a_i &= \frac{\tau_i + \sqrt{\tau_i^2 + 4\,\tau_iA_i\,(L_{i+1}-\mu)}}{2\,(L_{i+1}-\mu)}, \quad \tilde x_i = \frac{A_i\,y_i + a_i\,x_i}{A_i + a_i},
     \\
      y_{i+1} &=\arg\min_{u\in B_{\tau}^s}\left\{
        \ell_{\psi_s}(u;\tilde x_i) + \tfrac{L_{i+1}}2\|u-\tilde x_i\|^2
      \right\};
    \end{align}
    if 
    \begin{equation}
      \ell_{\psi_s}(y_{i+1};\tilde x_i) + \tfrac{(1-\chi)\,L_{i+1}}4\,\|y_{i+1}-\tilde x_i\|^2 
      \;\ge\;\psi_s(y_{i+1}),
    \end{equation}
    go to step~3; else set $L_{i+1}\leftarrow2\,L_{i+1}$ and repeat step~2;
  \item[{\bf 3.}] update
    \begin{align}
      A_{i+1} &= A_i + a_i, \qquad \tau_{i+1} = \tau_i + a_i\,\mu,\\
      s_{i+1} &= (L_{i+1}-\mu)\,(\tilde x_i - y_{i+1}),\\
      x_{i+1} &= \frac{1}{\tau_{i+1}}\left[\mu\,a_i\,y_{i+1} + \tau_i\,x_i - a_i\,s_{i+1}\right];
    \end{align}
  \item[{\bf 4.}] if 
    \begin{equation}
      \|y_{i+1}-x_0\|^2 \;\ge\;\chi\,A_{i+1}\,L_{i+1}\,\|y_{i+1}-\tilde x_i\|^2,
    \end{equation}
    then go to step~5; otherwise stop with {\bf failure};
  \item[{\bf 5.}] If the inequality
    \begin{equation}
      \|s_{i+1}\|\;\le\;\sigma\,\|y_{i+1}-x_0\|
    \end{equation}
holds, then stop with {\bf success} and output $(y,L)=(y_{i+1},L_{i+1})$; otherwise set $i\leftarrow i+1$ and return to step~1.
\end{itemize}
\noindent\rule[0.5ex]{1\columnwidth}{1pt}

Several remarks about ADAP-FISTA are now given. ADAP-FISTA 
is either able to successfully find a suitable solution of \eqref{psi s} or is unable to do so in at most 
\[
\mathcal O\left(\sqrt{L_{\psi_s}}\right)
\]
number of iterations, where $L_{\psi_s}$ is the Lipschitz constant of the gradient of the objective in \eqref{psi s}. If the objective in \eqref{psi s} is strongly convex, then ADAP-FISTA always succeeds in solving \eqref{psi s}. For more technical details on the type of solution that ADAP-FISTA aims to find, see Appendix B in \cite{monteiro2024low}.

\end{appendix}

\bibliographystyle{plain}

\bibliography{references}
\end{document}